\documentclass[reqno]{amsart}

\usepackage{amsmath,amssymb,pstricks,pst-node,calc}

\newlength{\mvadepth}

\newlength{\mvaheight}

\newlength{\mylen}

\newcommand{\ri}{\mathcal{O}}
\newcommand{\mi}{\mathcal{P}}
\newcommand{\mymatrix}[1]{%
 \settoheight{\mylen}{$\displaystyle \left[ \begin{array}{lll} #1 \end{array} \right] $}
 \addtolength{\mylen}{24pt}
 \rule[-0.48\mylen]{0pt}{\mylen} \displaystyle \left[ \begin{array}{lll} #1 \end{array} \right] }

\newcommand{\qtate}[2]{%
\settoheight{\mvaheight}{$\displaystyle #1 #2$}%
\settodepth{\mvadepth}{$\displaystyle #1 #2$}%
\displaystyle  #1 \left. \rule[-\mvadepth]{0pt}{\mvaheight} \right/ #2 }

\newcommand{\mycell}[1]{%
\settoheight{\mvaheight}{#1}%
\settodepth{\mvadepth}{#1}%
\addtolength{\mvaheight}{0.15cm}%
\addtolength{\mvadepth}{0.15cm}%
\rule[-\mvadepth]{0pt}{\mvaheight + \mvadepth} {#1}}

\newcommand{\eop}{\hspace{1cm} \rule{6pt}{6pt}}

\newtheorem{theorem}{Theorem}[section]
\newtheorem{lemma}[theorem]{Lemma}
\newtheorem{corollary}[theorem]{Corollary}
\newtheorem*{main theorem}{Main Theorem}

\numberwithin{equation}{section}

\begin{document}

\title{Affine pavings for affine Springer fibers for split elements in $PGL(3)$}
\author{Vincent Lucarelli}
\address{Department of Mathematics\\
         University of Chicago\\ 
         5734 South University Avenue\\
         Chicago, IL 60637} 
\email{vincent@math.uchicago.edu}
\date{September 5, 2003}
\subjclass[2000]{22E67; 14L40}

\begin{abstract}
This paper constructs pavings by affine spaces for the
affine Springer fibers for $PGL(3)$ obtained from regular compact elements
in a split maximal torus.  These pavings are constructed by intersecting
the affine Springer fiber with a non-standard paving of the affine
Grassmannian.
\end{abstract}

\maketitle
\section*{}

Let $k$ be a field, 
$F = k((\pi))$ be the field of 
formal Laurent series over $k$, $\ri = k[[\pi]]$ be
the subring of formal power series, and $\mi = \pi\ri$
be the maximal ideal.  Let $G = PGL_{3}(F)$ and $K = PGL_{3}(\ri)$.
We write $X$ for the affine Grassmannian $G/K$ and $A$ for the 
diagonal maximal torus in $G$.	

In section \ref{affine paving of X} we develop for each
integer $a \ge 0$ a non-standard paving of $X$ by affine spaces,
referred to as the $a$-paving of $X$.  These non-standard pavings are
constructed using the standard Iwahori
subgroup $I$ and a conjugate of it, $I^{a}$, that depends on $a$. 
Each affine space in an $a$-paving is a union of intersections
of $I$-orbits and $I^{a}$-orbits in $X$.  In particular, the affine spaces are
preserved by $I \cap I^{a}$ and hence by its subgroup $A(\ri)$. 
 
Each affine space in an $a$-paving contains exactly one 
element of $X_{*}(A)$ and thus these lattice points index the
affine spaces.  This is analogous to how affine spaces are indexed 
in the standard paving of $X$ by $I$-orbits.  
In fact, when $a=0$ the $a$-paving of $X$
is identical to the standard affine
paving of $X$ by $I$-orbits. 

An $a$-paving differs from the standard paving of $X$ by $I$-orbits
in how the closure of each affine space relates to the paving.  
In the standard paving, the closure of an 
affine space is the union of smaller dimensional affine spaces.
For all $a$ other than zero, our $a$-paving is a paving in a weaker 
sense.  The closure of each affine space is 
not necessarily the union of other affine spaces in the paving.  However,
we can order the affine spaces $\mathbb{A}_{0}, \mathbb{A}_{1}, \mathbb{A}_{2}, \ldots$
so that $\mathbb{A}_{0} \cup \dots \cup \mathbb{A}_{n}$ is closed 
for all $n$.

Our non-standard pavings are interesting because they induce
affine pavings of certain fixed point sets in $X$.  Specifically,
when $\gamma$ is a regular element in $A(\ri)$ the set
\[ X^{\gamma} = \{ g \in G/K : \gamma g = g \}\]
admits a paving by affine spaces when intersected  
with a particular $a$-paving of $X$ that is determined by $\gamma$.
The set $X^{\gamma}$ is called an affine Springer fiber.
They were first studied by Kazhdan and Lusztig in \cite{KL88}.

The value of $a$ that induces an affine paving of $X^{\gamma}$
is determined as follows.  We use $v(x)$ to denote the valuation of 
$x \in F$ and we take $v(0) = +\infty$.
Now $\gamma$ has the form
\[ \gamma = \mymatrix{ u_{1} & & \\ & u_{2} & \\ & & u_{3}} \]
where $u_{1}$, $u_{2}$, and $u_{3}$ are in $\ri^{\times}$ and distinct.
We can permute the entries of $\gamma$ so that 
\[ v\left( 1 - \frac{u_{1}}{u_{2}}\right) = v\left( 1 - \frac{u_{1}}{u_{3}}\right) = m \]
and
\[ v\left( 1 - \frac{u_{2}}{u_{3}}\right) = n \]
with $n \ge m \ge 0$.  Then $a = n - m$.

\begin{main theorem}
The intersection of $X^{\gamma}$ with the $a$-paving of $X$ determined 
by $a = n-m$ yields an affine paving of $X^{\gamma}$.  In particular,
the intersection of $X^{\gamma}$ with each affine space in the $a$-paving
of $X$ is an affine space.
\end{main theorem}

To prove the result,
we explicitly calculate the intersection of $X^{\gamma}$ with the 
affine spaces of the $a$-paving of $X$ in section \ref{fixed point calculations}.
The precise statement of the main theorem is Theorem \ref{precise main theorem}.

Many other authors have found affine pavings of affine Springer fibers
for specific groups with an equivalued condition on the valuation of roots \cite{Fan96} \cite{LW}
\cite{LS91} \cite{Sag00} \cite{Som97}.
Goresky, Kottwitz, and MacPherson
proved a general result that gives an affine paving of affine Springer fibers for any 
connected reductive group assuming the equivalued condition \cite{GKM}.  This paper develops the 
first known affine paving of an affine Springer fiber in the non-equivalued case.

We conjecture that this method can be extended to prove the same result
when $k$ is algebraically closed, $\mathop{char}(k) \not= 2$, and the maximal
torus splits as $E^{\times} \times F^{\times}$ where $E/F$ is a quadratic
extension.  Provided this is true, in $PGL(3)$ an affine paving of $X^{\gamma}$ is known 
 for all regular semisimple $\gamma$ 
except in characteristic 2 and 3.
In the case of the torus $E^{\times}$ where $E/F$ is a cubic
extension, $n$ and $m$ are forced to be equal  and 
in characteristic other than 2 and 3 \cite{GKM} applies.
(\cite{GKM} actually applies to the Lie algebra of $G$, but the result
is equivalent for the group.)

I wish to thank my advisor Robert Kottwitz for suggesting this problem
and for his generous help.

\section{Notation, Definitions, \& Precise Statement of the Main Theorem}

An element of $X$ that has a diagonal matrix coset representative 
can be expressed in terms of a diagonal matrix with monomial entries.
That representative is equivalent in $G$ to an element of the form
\[ \mymatrix{ 1 & & \\ & \pi^s & \\ & & \pi^t}. \]
We denote such elements of $X$ by the coordinates $(s,t)$.  
(These points are the vertices in the main apartment of the building for $G$.) 

Let $I$ denote the standard Iwahori subgroup
\[ I = \mymatrix{ \ri^{\times} & \ri & \ri \\ \mi & \ri^{\times} & \ri \\ \mi & \mi & \ri^{\times}} \]
and for $a \ge 0$ define the conjugate
\begin{eqnarray*} 
I^{a} & = & \mymatrix{ 1 & & \\ & \pi^{a} & \\ 
& & \pi^{a}} I \mymatrix{ 1 & & \\ & \pi^{-a} & \\ & & \pi^{-a}} \\
& = & 
\mymatrix{ \ri^{\times} & \mi^{-a} & \mi^{-a} \\ 
\mi^{a+1} & \ri^{\times} & \ri \\ \mi^{a+1} & \mi & \ri^{\times}}.
\end{eqnarray*}

\begin{figure}
	\begin{center}
    \psset{xunit=0.166666666666667 in}
    \psset{yunit=0.288675134594813 in}
    {\scriptsize  
    \begin{pspicture}(-1,-1)(13,7)
	
	\pspolygon[fillstyle=solid,fillcolor=gray,linestyle=none](6,2)(5,3)(4,2)
	\uput{0.4}[150](6,2){$I^{0}$}

	\psline[linewidth=0.5pt](2.25,6)(2.25,6.25)(9.75,6.25)(9.75,6)
	\rput[b]{0}(6,6.35){type $1^{0}$}
	\rput[l]{60}(10,6){\ type $2^{0}$}
	\psline(10.5,6)(11,6.25)(12.25,6.25)(12.25,2.25)(12,2.25)
	\rput[l]{0}(12.25,6.25){\ type $3^{0}$}
	\rput[l]{0}(12,2){\ type $4^{0}$}
	\psline(12,1.75)(12.25,1.75)(12.25,-0.25)(9,-0.25)(8.5,0)
	\rput[l]{0}(12.25,-0.25){\ type $5^{0}$}
	\rput[l]{-60}(8,0){\ type $6^{0}$}
	\psline[linewidth=0.5pt](4.25,0)(4.25,-0.25)(7.75,-0.25)(7.75,0)
	\rput[t]{0}(6,-0.35){type $7^{0}$}	
	\rput[r]{60}(4,0){type $8^{0}$\ }
	\psline[linewidth=0.5pt](3.5,0)(3,-0.25)(-0.25,-0.25)(-0.25,1.75)(0,1.75)
	\rput[r]{0}(-0.25,-0.25){type $9^{0}$\ }
	\rput[r]{0}(0,2){type $10^{0}$\ }
	\psline[linewidth=0.5pt](0,2.25)(-0.25,2.25)(-0.25,6.25)(1,6.25)(1.5,6)
	\rput[r]{0}(-0.25,6.25){type $11^{0}$\ }
	\rput[r]{-60}(2,6){type $12^{0}$\ }

	\psline[linewidth=0.1pt](0,2)(2,0)
	\psline[linewidth=0.1pt](0,4)(4,0)
	\psline[linewidth=0.1pt](0,6)(6,0)
	\psline[linewidth=1.5pt](2,6)(8,0)
	\psline[linewidth=0.1pt](4,6)(10,0)
	\psline[linewidth=0.1pt](6,6)(12,0)
	\psline[linewidth=0.1pt](8,6)(12,2)
	\psline[linewidth=0.1pt](10,6)(12,4)
	\psline[linewidth=0.1pt](12,6)(12,6)
	\psline[linewidth=0.1pt](0,4)(2,6)
	\psline[linewidth=0.1pt](0,2)(4,6)
	\psline[linewidth=0.1pt](0,0)(6,6)
	\psline[linewidth=0.1pt](2,0)(8,6)
	\psline[linewidth=1.5pt](4,0)(10,6)
	\psline[linewidth=0.1pt](6,0)(12,6)
	\psline[linewidth=0.1pt](8,0)(12,4)
	\psline[linewidth=0.1pt](10,0)(12,2)
	\psline[linewidth=0.1pt](12,0)(12,0)
	\psline[linewidth=0.1pt](0,0)(12,0)
	\psline[linewidth=0.1pt](0,1)(12,1)
	\psline[linewidth=1.5pt](0,2)(12,2)
	\psline[linewidth=0.1pt](0,3)(12,3)
	\psline[linewidth=0.1pt](0,4)(12,4)
	\psline[linewidth=0.1pt](0,5)(12,5)
	\psline[linewidth=0.1pt](0,6)(12,6)
	
	\psdots[dotstyle=*](0,0) 
	\psdots[dotstyle=*](0,2) \uput[330](0,2){$(3,3)$}
	\psdots[dotstyle=*](0,4) \uput[330](0,4){$(2,4)$}
	\psdots[dotstyle=*](0,6) \uput[330](0,6){$(1,5)$}
	\psdots[dotstyle=*](3,1) \uput[330](3,1){$(2,1)$}
	\psdots[dotstyle=*](3,3) \uput[330](3,3){$(1,2)$}
	\psdots[dotstyle=*](3,5) \uput[330](3,5){$(0,3)$}
	\psdots[dotstyle=*](6,0) 
	\psdots[dotstyle=*](6,2) \uput[330](6,2){$(0,0)$}
	\psdots[dotstyle=*](6,4) \uput[330](6,4){$(-1,1)$}
	\psdots[dotstyle=*](6,6) \uput[330](6,6){$(-2,2)$}
	\psdots[dotstyle=*](9,1) \uput[330](9,1){$(-1,-2)$}
	\psdots[dotstyle=*](9,3) \uput[330](9,3){$(-2,-1)$}
	\psdots[dotstyle=*](9,5) \uput[330](9,5){$(-3,0)$}
	\psdots[dotstyle=*](12,0) 
	\psdots[dotstyle=*](12,2) 
	\psdots[dotstyle=*](12,4) 
	\psdots[dotstyle=*](12,6) 
	\end{pspicture}}
    \end{center}
\caption{Example for $a = 0$}\label{I stab}
\end{figure}
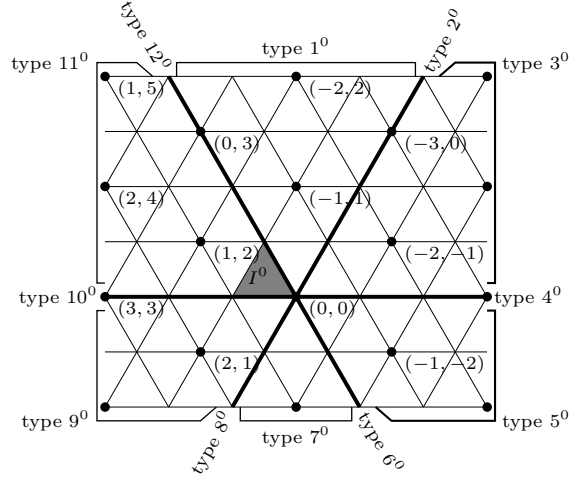

The subgroup $I^{a}$ acts on points of $X$ by left multiplication.  This action fixes
the point $(a,a)$, which we call the base point relative to $a$.
For a particular $I^{a}$,
the other points $(s,t)$ in the main apartment of $X$ can be divided into 
12 types based on the relationship of $a$ and the coordinates $(s,t)$.
\begin{center}
\begin{tabular}{|l|l|}\hline
Type & Condition \\ \hline
\mycell{$1^{a}$} & \mycell{ $s < a < t$ } \\ \hline
\mycell{$2^{a}$} & \mycell{ $s < t = a$ } \\ \hline
\mycell{$3^{a}$} & \mycell{ $s < t < a$ } \\ \hline
\mycell{$4^{a}$} & \mycell{ $s = t < a$ } \\ \hline
\end{tabular}
\hspace{1cm}
\begin{tabular}{|l|l|}\hline
Type & Condition \\ \hline
\mycell{$5^{a}$} & \mycell{ $t < s < a$ } \\ \hline
\mycell{$6^{a}$} & \mycell{ $t < s = a$ } \\ \hline
\mycell{$7^{a}$}  & \mycell{ $t < a < s$ } \\ \hline
\mycell{$8^{a}$}  & \mycell{ $a = t < s$ } \\ \hline
\end{tabular}
\hspace{1cm}
\begin{tabular}{|l|l|}\hline
Type & Condition \\ \hline
\mycell{$9^{a}$}  & \mycell{ $a < t < s$ } \\ \hline
\mycell{$10^{a}$} & \mycell{ $a < t = s$ } \\ \hline
\mycell{$11^{a}$} & \mycell{ $a < s < t$ } \\ \hline
\mycell{$12^{a}$} & \mycell{ $a = s < t$ } \\ \hline
\end{tabular}
\end{center}
Figure \ref{I stab} shows how the vertices in the main
apartment are partitioned into the twelve types relative to
$I^{0}$.  For simplicity, we only label
the vertices corresponding to elements in $SL_{3}(F)$.

Every point in $X$ is in the $I^{a}$-orbit of some point in the main apartment.
This gives a decomposition of $X$ into disjoint sets
\[ X = \bigsqcup_{ x \in \textrm{vert}} \qtate{ I^{a} x K}{K} \]
where $\textrm{vert}$ is the set of vertices $(s,t)$. 

To pave $X$ by affine spaces, we first divide $X$ into three disjoint sets
\[ S = \bigsqcup_{\textrm{$x$ type $1^0$, $2^0$, or $3^0$}} \qtate{IxK}{K}, \hspace{1in}
 T = \bigsqcup_{\textrm{$x$ type $5^0$, $6^0$, or $7^0$}} \qtate{IxK}{K}, \]
and
\[ V = X \setminus (S \cup T).\]
The intersection of the $I^{a}$-orbit of a vertex $v$ with the 
set $S$ 
is denoted by 
\[ S^{a}_{v} = S \cap \qtate{I^{a} v K}{K}\]
with analogous notation $T^{a}_{v}$ for the set $T$
and $V^{a}_{v}$ for the set V.
In the course of proving the main theorem, we will show that 
the sets $V^{0}_{v}$, $S^{a}_{v}$ and $T^{a}_{v}$ are affine spaces
and together they give a non-standard paving of $X$ by affine spaces.

We now state the main theorem.  
\begin{theorem}\label{precise main theorem}
Let $a=n-m$.  
The sets $X^{\gamma}\cap V_{v}^{0}$, $X^{\gamma}\cap S_{v}^{a}$, $X^{\gamma}\cap T_{v}^{a}$
are affine spaces that 
form an affine paving of $X^{\gamma}$ as $v$ ranges over the vertices in the 
main apartment of $X$.  
\end{theorem}

\section{Understanding $I^{a}$-orbits}

The sets $S$ and $T$ are defined by $I$-orbits of vertices,
but the main theorem utilizes sets that are defined in part by the $I^{a}$-orbits of those vertices.  
This motivates us to analyze the relationship of $I$-orbits and $I^{a}$-orbits.

To begin, we find a unique matrix coset representative for each point 
in the $I^{a}$-orbit of a vertex $(s,t)$.  With this enumeration of points, 
we can explicitly describe how $I$ and $I^{a}$-orbits compare.

\subsection{Enumerating Cosets}

Consider the $I^{a}$-orbit of the vertex $x=(s,t)$,
\[ \qtate{ I^{a} x K}{K} = \qtate{I^{a}}{I^{a} \cap xKx^{-1}}.\]
For an arbitrary matrix in $k \in K$ the conjugate $xkx^{-1}$ is 
\[ x\mymatrix{a & b & c \\d & e & f \\g & h & i } x^{-1} = 
  \mymatrix{a & b\pi^{-s} & c\pi^{-t} \\d\pi^{s} & e & f\pi^{s-t} \\g\pi^{t} & h\pi^{t-s} & i}. \]
We are interested in the intersection 
\[ I^{a}_{(s,t)} = I^{a} \cap xKx^{-1}. \]
The matrix form of this intersection depends upon the type of $x$ 
relative to $a$.  For example, when $x$ is type $1^{a}$ ($s < a < t$) 
\[ I^{a}_{(s,t)} = \mymatrix{\ri^{\times} & \mi^{-s} & \mi^{-a} \\ 
	\mi^{a+1} & \ri^{\times} & \ri \\ \mi^{t} & \mi^{t-s} & \ri^{\times} }\]
and
\begin{eqnarray*}
&& \qtate{I^{a}}{I^{a}_{(s,t)}}  =  
\qtate{ \mymatrix{ \ri^{\times} & \mi^{-a} & \mi^{-a} \\ \mi^{a+1} & \ri^{\times} & \ri \\ \mi^{a+1} & \mi & \ri^{\times}}}{
	\mymatrix{\ri^{\times} & \mi^{-s} & \mi^{-a} \\ \mi^{a+1} & \ri^{\times} & \ri \\ \mi^{t} & \mi^{t-s} & \ri^{\times} }}
	\\
& = & \qtate{ \mymatrix{1 & \mi^{-a} & 0 \\ 0 & 1 & 0 \\ \mi^{a+1} & \mi & 1}\mymatrix{\ri^{\times} & \mi^{-s} & \mi^{-a} \\ \mi^{a+1} & \ri^{\times} & \ri \\ \mi^{t} & \mi^{t-s} & \ri^{\times} }}{\mymatrix{\ri^{\times} & \mi^{-s} & \mi^{-a} \\ \mi^{a+1} & \ri^{\times} & \ri \\ \mi^{t} & \mi^{t-s} & \ri^{\times} }} \\
& = & \qtate{ \mymatrix{1 & \mi^{-a} & 0 \\ 0 & 1 & 0 \\ \mi^{a+1} & \mi & 1}}{\mymatrix{1 & \mi^{-s} & 0 \\ 0 & 1 & 0 \\ \mi^{t} & \mi^{t-s} & 1}}.
\end{eqnarray*}
Then for each $M \in I^{a}$, we can factor $M$ as 
\[  \mymatrix{ 1 & i & 0 \\ 0 & 1 & 0 \\ y & z & 1} \textrm{ times a matrix in} 
	\mymatrix{\ri^{\times} & \mi^{-s} & \mi^{-a} \\ \mi^{a+1} & \ri^{\times} & \ri \\ \mi^{t} & \mi^{t-s} & \ri^{\times} }
	\]
where 
\begin{eqnarray*}
	i = i_{-a}\pi^{-a} + \dots + i_{-s-1}\pi^{-s-1} \\
	y = y_{a+1}\pi^{a+1} + \dots + y_{t-1}\pi^{t-1} \\
	z = z_{1}\pi^{1} + \dots + z_{t-s-1}\pi^{t-s-1}.
\end{eqnarray*}
We indicate $i$, $y$, and $z$ are such polynomials by writing
$i \in \qtate{\mi^{-a}}{\mi^{-s}}$, $y \in \qtate{\mi^{a+1}}{\mi^{t}}$, 
and $z \in \qtate{\mi}{\mi^{t-s}}$.  

So each coset in $\qtate{I^{a}}{I^{a}_{(s,t)}}$ (when $s < a < t$) can be represented using a triple $(i,y,z)$
by factoring $M$ as above.  In fact, each triple represents a unique coset.  To see this, 
suppose the cosets represented by $(i,y,z)$ and $(i',y',z')$ are equivalent, then
\[ \mymatrix{ 1 & -i & 0 \\ 0 & 1 & 0 \\ -y & iy - z & 1} 
	\mymatrix{ 1 & i' & 0 \\ 0 & 1 & 0 \\ y' & z' & 1} 
	= \mymatrix{ 1 & i' - i & 0 \\ 0 & 1 & 0 \\ y' - y & y(i - i') + (z' - z) & 1}
	\in I^{a}_{(s,t)}.\]
This forces $i = i'$ because $i' - i \in \mi^{-s}$ but the valuations of 
$i$ and $i'$ are less than $-s$.  Similarly, we have $y = y'$ and $z' = z$.
So the set of matrices of the form
\[  \mymatrix{ 1 & i & 0 \\ 0 & 1 & 0 \\ y & z & 1} \]
with $i \in \qtate{\mi^{-a}}{\mi^{-s}}$, $y \in \qtate{\mi^{a+1}}{\mi^{t}}$, 
and $z \in \qtate{\mi}{\mi^{t-s}}$ gives a complete 
set of unique coset representatives for the elements in $\qtate{I^{a}}{I^{a}_{(s,t)}}$ when $s < a < t$. 
This description also gives a complete
enumeration of all the points in the $I^{a}$-orbit of a vertex of type $1^{a}$.   

A similar analysis can be performed for each vertex type.  Using 
\[ \mymatrix{ 1 & i & j \\ x & 1 & k \\ y & z & 1} \] as a 
generic coset representative, the following table summarizes
how to enumerate, in a standard way,
the points in the $I^{a}$-orbit of each vertex type.  
\begin{center}
{\small
\begin{tabular}{|l|c|c|c|c|c|c|}\hline
Type &  $i$ & $j$ & $k$ & $x$ & $y$ & $z$ \\ \hline
\mycell{ $1^a$} & \mycell{$\qtate{\mi^{-a}}{\mi^{-s}}$} & 0 & 0 & 0 & 
		\mycell{$\qtate{\mi^{a+1}}{\mi^{t}}$} & \mycell{ $\qtate{\mi}{\mi^{t-s}}$} \\ \hline
\mycell{ $2^a$} & \mycell{$\qtate{\mi^{-a}}{\mi^{-s}}$} & 0 & 0 & 0 & 0 & 
		\mycell{ $\qtate{\mi}{\mi^{t-s}}$} \\ \hline
\mycell{ $3^a$} & \mycell{$\qtate{\mi^{-a}}{\mi^{-s}}$} & 
		\mycell{$\qtate{\mi^{-a}}{\mi^{-t}}$} & 0 & 0 & 0 & \mycell{ $\qtate{\mi}{\mi^{t-s}}$}\\ \hline
\mycell{ $4^a$} & \mycell{$\qtate{\mi^{-a}}{\mi^{-s}}$} & 
		\mycell{$\qtate{\mi^{-a}}{\mi^{-t}}$} & 0 & 0 & 0 & 0 \\ \hline
\mycell{ $5^a$} & \mycell{$\qtate{\mi^{-a}}{\mi^{-s}}$} & 
		\mycell{$\qtate{\mi^{-a}}{\mi^{-t}}$} & \mycell{$\qtate{\ri}{\mi^{s-t}}$} & 0 & 0 & 0\\ \hline
\mycell{ $6^a$} & 0 & \mycell{$\qtate{\mi^{-a}}{\mi^{-t}}$} & 		\mycell{$\qtate{\ri}{\mi^{s-t}}$} & 0 & 0 & 0 \\ \hline
\mycell{ $7^a$} & 0 & \mycell{$\qtate{\mi^{-a}}{\mi^{-t}}$} & 		\mycell{$\qtate{\ri}{\mi^{s-t}}$} & \mycell{$\qtate{\mi^{a+1}}{\mi^{s}}$} & 0 & 0\\ \hline
\mycell{ $8^a$} & 0 & 0 & \mycell{$\qtate{\ri}{\mi^{s-t}}$} & 
		\mycell{$\qtate{\mi^{a+1}}{\mi^{s}}$} & 0 & 0\\ \hline
\mycell{ $9^a$} & 0 & 0 & \mycell{$\qtate{\ri}{\mi^{s-t}}$} & 
		\mycell{$\qtate{\mi^{a+1}}{\mi^{s}}$} & \mycell{$\qtate{\mi^{a+1}}{\mi^{t}}$} & 0\\ \hline
\mycell{ $10^a$} & 0 & 0 & 0 & \mycell{$\qtate{\mi^{a+1}}{\mi^{s}}$} & 
		\mycell{$\qtate{\mi^{a+1}}{\mi^{t}}$} & 0\\ \hline
\mycell{ $11^a$} & 0 & 0 & 0 & \mycell{$\qtate{\mi^{a+1}}{\mi^{s}}$} & 
		\mycell{$\qtate{\mi^{a+1}}{\mi^{t}}$} & \mycell{ $\qtate{\mi}{\mi^{t-s}}$} \\ \hline
\mycell{ $12^a$} & 0 & 0 & 0 & 0 & \mycell{$\qtate{\mi^{a+1}}{\mi^{t}}$} & 
		\mycell{ $\qtate{\mi}{\mi^{t-s}}$} \\ \hline
\end{tabular}}
\end{center}
A matrix $M$ that represents a point in the $I^{a}$-orbit of a vertex of type $r$
is said to be in standard form if the diagonal entries of $M$ are one and the 
other entries of $M$ follow the specification for type $r$ in the above table.
When referring to a point $M$, we will also refer to particular values
in the standard form of $M$ via the variables
$i_{M},j_{M},k_{M},x_{M},y_{M},z_{M}$.  For example, if $M$ is a 
point in the $I^{a}$-orbit of a vertex of type $1^{a}$, then the 
$j_{M} = 0$ and $z_{M} \in \qtate{\mi}{\mi^{t-s}}$.

\subsection{Stationary \& Non-stationary Points}\label{affine paving of X}

Given a point $M$ in the $I$-orbit of a vertex $v$, we want to find the $I^{a}$-orbit that
contains $M$.  If $M$ is in the $I^{a}$-orbit of $v$, we call
the point stationary; otherwise the point is called non-stationary.  When $a = 0$,
all points are stationary and the $a$-paving is identical to the paving
of $X$ by $I$-orbits.  So throughout this section we assume $a > 0$.   

Since we only consider the intersection of $I^{a}$-orbits with the sets $S$ and $T$,
we only need to compare the $I$ and $I^{a}$-orbits for vertices of 
type $1^{0}$, $2^{0}$, $3^{0}$, $5^{0}$, $6^{0}$, and $7^{0}$.

\subsubsection{$v$ is Type $1^{0}$, $2^{0}$, or $3^{0}$}\label{alt retract I}

This section is devoted to proving the following lemma.
Recall that we denote the valuation
of the $y$ component of the standard enumeration of the point $M$ in the $I$-orbit
of a vertex by $v(y_{M})$.

\begin{lemma}\label{123 movement}
Let $M$ be a point in the $I$-orbit of the vertex $v = (s,t)$.
\begin{itemize}
\item[(a)] If $v$ is type $2^{0}$ or $3^{0}$ then $M$ is stationary
\item[(b)] If $v$ is type $1^{0}$ and $v( y_{M}) > a$ then $M$ is stationary
\item[(c)] If $v$ is type $1^{0}$ and $v( y_{M}) \le a$ then $M$ is in the $I^{a}$-orbit
	of $w = (s - d, t-2d)$ where $d = t - v(y_{M})$
\end{itemize}
\end{lemma}

\emph{Part (a) \& (b)}: $M$ is stationary if $M \in I \cap I^{a}$.
Since 
\begin{eqnarray*}
I \cap I^{a} & = & \mymatrix{ \ri^{\times} & \ri & \ri \\ 
\mi & \ri^{\times} & \ri \\ \mi & \mi & \ri^{\times}}
\cap
\mymatrix{ \ri^{\times} & \mi^{-a} & \mi^{-a} \\ 
\mi^{a+1} & \ri^{\times} & \ri \\ \mi^{a+1} & \mi & \ri^{\times}}  \\
& = & \mymatrix{ \ri^{\times} & \ri & \ri \\ 
\mi^{a+1} & \ri^{\times} & \ri \\ \mi^{a+1} & \mi & \ri^{\times}}
\end{eqnarray*}
we see that any $M$ representing a point in the $I$-orbit
of a vertex of type $2^{0}$, $3^{0}$, or $1^{0}$ with
$v(y_M) > a$ is in $I \cap I^{a}$ and so
$M$ is stationary.

\emph{Part (c)}:\   Since $0 < v(y_{M}) < t$, we have $d > 0$,
which can be thought of as the distance between $v$ and $w$.
To show that  
$M$ is in the $I^{a}$-orbit of $w$, it suffices to find $M'$ in $I^{a}$ such
that 
\[ M' w \in M v K.\]
For simplicity we wish to work with matrices
rather than matrices modulo scalars.
Since $\det(w) = s + t - 3d$ and $\det(v) = s + t$,
we replace $w$ with the equivalent element $\pi^{d}w$
and thus consider
\[ \pi^{d}M' w \in M v K.\]
Rearranging terms yields
\[ M^{-1} M' \in \pi^{-d}vKw^{-1}. \]
Since $\pi^{-d}vKw^{-1}$ is the set of matrices of the form
\[ \mymatrix{ \mi^{-d} &  \mi^{-s} & \mi^{d-t} \\ \mi^{s-d} &  \ri & \mi^{s-t+d} \\  
	\mi^{t-d} &  \mi^{t-s} & \mi^{d} } \]
with determinant a unit and $\det( M^{-1} M')$ is a unit,
the condition on $M'$ reduces to 
\[ M^{-1} M' \in \mymatrix{ \mi^{-d} &  \mi^{-s} & \mi^{d-t} \\ \mi^{s-d} &  \ri & \mi^{s-t+d} \\  
\mi^{t-d} &  \mi^{t-s} & \mi^{d} }.\]
To compute the matrix form of the left side, we need to know 
the standard form for $M'$, which is determined by $w$'s type relative
to $I^{a}$.  Later, it will be convenient to also know $w$'s type
relative to $I$.  

\begin{lemma}\label{w is type 3a}
Suppose that $v=(s,t)$ is type $1^0$ and 
$1 \le v(y) \le \min(a, t-1)$.  Define $d = t - v(y)$.  Then 
$ w = (s - d, t - 2d) $ is type $3^a$ and either type $1^0$, $2^0$, or $3^0$.
\end{lemma}

\emph{Proof}.\  The inequalities
\[ (s - d) < 0 \quad\textrm{and}\quad (s-d) < (t - 2d)\]
imply $w$ is type $1^0$, $2^0$, or $3^0$
and
\[ (s - d) < (t - 2d) < a\]
implies $w$ is type $3^{a}$.

By assumption, $v(y) \le \min(a, t-1)$
so $v(y) < t$ and thus $d \ge 1$.  Since $v$ is 
type $1^0$ we know $s < 0$.  Therefore the inequality $(s - d) < 0$ is clear.

The inequality $(s - d) < (t - 2d)$ is equivalent to $s < v(y)$
which follows from $s < 0 < 1 \le v(y)$.
Finally, the inequality $(t - 2d) < a$ is equivalent to $2v(y) < a + t$
which follows from $v(y) \le a$ and $v(y) \le t - 1 < t$.~\eop

We can now compute $M^{-1}M'$.  To simplify notation,
we use $i = i_{M}$, $i' = i_{M'}$, and similar simplifications
for the other variables.  Note that 
\[ M = \mymatrix{ 1 & i & 0 \\ 0 & 1 & 0 \\ y & z & 1}.\]

Because $w$ is type $3^{a}$ we will look for $M' \in I^{a}$ of the form
\[ M' =  \mymatrix{1 & i' & j' \\ 0 & 1 & 0 \\ 0 & z' & 1}. \]
Then
\begin{eqnarray*}
 M^{-1}M' & = & \mymatrix{ 1 & -i & 0 \\ 0 & 1 & 0 \\ -y & iy - z & 1 } \mymatrix{1 & i' & j' \\ 0 & 1 & 0 \\ 0 & z' & 1} \\
	& = & \mymatrix{ 1 & (i' - i) & j' \\  	0 & 1 & 0 \\	-y & [ y(i-i') + (z' - z)] & (1 - yj') }.
\end{eqnarray*}
Therefore, we need
\begin{equation}\label{m hat m condition I}
 \mymatrix{ 1 & (i' - i) & j' \\  	0 & 1 & 0 \\	-y & [ y(i-i') + (z' - z)] & (1 - yj') } 
 \subset \mymatrix{ \mi^{-d} &  \mi^{-s} & \mi^{d-t} \\ \mi^{s-d} &  \ri & \mi^{s-t+d} \\  \mi^{t-d} &  \mi^{t-s} & \mi^{d} } \end{equation}
for $M'$ to represent the point $M$ in the $I^{a}$-orbit of $w$.

From our hypotheses for part (c), we have
\[ i \in \qtate{\ri}{\mi^{-s}}, \quad y \in \qtate{\mi}{\mi^{t}}, \textrm{\quad and \quad} 
 z \in \qtate{\mi}{\mi^{t-s}} \]
with $v(y) \le \min(a, t-1)$.  To determine $M'$, first 
express $z$ as a polynomial in $\pi$
\[ z = z_{1} \pi^{1} + z_{2} \pi^{2} + \dots + z_{t-s-1} \pi^{t-s-1}. \]
and define
\[ \lceil z \rceil = z_{t-s-d} \pi^{t-s-d} + \dots + z_{t-s-1} \pi^{t-s-1}. \]
Let
\[\begin{array}{rclcr}  z' & = & z - \lceil z \rceil &\in& \qtate{\mi}{\mi^{t-s-d}}\\  
i' & = &\displaystyle i - \frac{\lceil z \rceil}{y}  &\in& \qtate{\ri}{\mi^{d-s}}\\  
j' & = &\displaystyle \frac{1}{y}  &\in& \qtate{\mi^{d-t}}{\mi^{2d-t}}.
\end{array}\]

To verify that $ M^{-1}M' \in \pi^{-d}vKw^{-1}$, we check the 
interesting matrix entries in equation (\ref{m hat m condition I}). 
\begin{center}
\begin{tabular}{|c|c|}\hline
Entry & Check \\ \hline
$(1,2)$ & \mycell{$\displaystyle v(i' - i) = v\left( \frac{\lceil z \rceil}{y}\right) = t-s-d-(t-d) = -s$ so 
	$(i' - i) \in \mi^{-s}$} \\ \hline
$(1,3)$ & \mycell{$\displaystyle v(j') = -v(y) = d-t \textrm{ so } j' \in \mi^{d-t}$} \\ \hline
$(3,1)$ & \mycell{$\displaystyle v(y) = t-d \textrm{ so } -y \in \mi^{t-d}$} \\ \hline
$(3,2)$ & \mycell{$ \begin{array}{lcl} \displaystyle y(i - i') + (z' - z) 
    & = & \displaystyle y\left( \frac{\lceil z \rceil}{y} \right) - \lceil z \rceil \\    
    & = & 0\end{array}$} \\ \hline
$(3,3)$ & \mycell{$\begin{array}{lcl}\displaystyle  (1 - yj') 
    & = & \displaystyle  1 -  y\left( \frac{1}{y} \right) \\  
    & = & 0\end{array}$} \\ \hline
\end{tabular}
\end{center}

Thus the point $M$ in the $I$-orbit of $v$ is also represented by $M'$ in
the $I^{a}$-orbit of $w$.~\eop

We noted that for a non-stationary 
point the corresponding $w$ was of type $1^{0}$, $2^{0}$, or $3^{0}$.  
This shows that when we rearrange elements of $S$ by $I^{a}$-orbits, any 
non-stationary point is in the $I^{a}$-orbit of a vertex in $S$.

\begin{corollary}
\[ S = \bigsqcup_{\textrm{$v$ type $1^0$, $2^0$, or $3^0$}} S^{a}_{v} \]
\end{corollary}

\emph{Proof}.  It is clear from the definition of $S^{a}_{v}$
that $S^{a}_{v} \subset S$.

Let $M \in S$.  Then for some $v$ of type $1^0$, $2^0$, or $3^0$
\[ M \in \qtate{IvK}{K}.\]

If $M$ is stationary, then $M \in S^{a}_{v}$.  Otherwise, $M$ 
is non-stationary and there is a $w$ of type $1^0$, $2^0$, or $3^0$
such that $M \in S^{a}_{w}$.~\eop

\subsubsection{Structure of $S^{a}_{v}$ for $v$ of Type $1^{0}$, $2^{0}$, or $3^{0}$ }\label{alt retract II}

In the last section, we categorized the points of $S$ as stationary
and non-stationary.  Given a vertex $v$ of type $1^{0}$, $2^{0}$, or $3^{0}$,
we already know which points in $S^{a}_{v}$ are stationary.  To complete
the description of $S^{a}_{v}$, we need to determine the 
non-stationary points in the set.

When $v$ is type $1^{a}$ or $2^{a}$ the only points in $S^{a}_{v}$
are stationary points because non-stationary move to $I^{a}$-orbits of type $3^{a}$ vertices.  
Since only type $1^{0}$ vertices can be type $1^{a}$ or $2^{a}$,
we can simply apply the stationary condition
$v(y) > a$ to the elements in the enumeration of points in the $I$-orbit of a
vertex of type $1^{0}$ to find $S^{a}_{v}$ in this case.  Therefore, 
when $v$ is type $1^{a}$ or $2^{a}$, we have 
$s < a \le t$ and $S^{a}_{v}$ is
enumerated by
\begin{eqnarray*}
	i & \in & \qtate{\ri}{\mi^{-s}} \\
	y & \in & \qtate{\mi^{a+1}}{\mi^{t}} \\
	z & \in & \qtate{\mi}{\mi^{t-s}}. 
\end{eqnarray*}

When $v$ is type $3^{a}$, the situation is more complicated
because non-stationary points move into $I^{a}$-orbits of these vertices.  Further 
complicating the situation is the fact that $v$ can be 
type $1^{0}$, $2^{0}$, or $3^{0}$ and so the 
enumeration of the stationary points is different in each case.

To list the stationary points, first suppose that $v$ is 
type $3^{a}$ and either type $1^{0}$ or type $2^{0}$.
Then $s < 0 \le t < a$.  If $v$ is type $2^0$, then $y = 0$.
Otherwise the stationary condition $v(y) > a$
applied to $y \in \qtate{\mi}{\mi^{t}}$ forces $y = 0$.  So 
the stationary points in $S^{a}_{v}$ are given by
\begin{eqnarray*}
	i & \in & \qtate{\ri}{\mi^{-s}} \\
	z & \in & \qtate{\mi}{\mi^{t-s}}. 
\end{eqnarray*}

When $v$ is type $3^{a}$ and $3^{0}$, 
$s < t < 0 < a$ and all the points are
stationary
\begin{eqnarray*}
	i & \in & \qtate{\ri}{\mi^{-s}} \\
	j & \in & \qtate{\ri}{\mi^{-t}} \\
	z & \in & \qtate{\mi}{\mi^{t-s}}. 
\end{eqnarray*}

For each positive integer $d$, there are non-stationary
points in $S$ that contribute to the $I^{a}$-orbit 
of $v = (s,t)$ if $v$ is type $3^{a}$
and
\begin{itemize}
\item[(a)] $w = (s + d, t + 2d)$ is a type $1^{0}$
\item[(b)] there exists $M \in \qtate{IwK}{K}$ so that $v(y) \le a$ and $d = (t + 2d) - v(y)$ (where $y \equiv y_{M}$).
\end{itemize}
Condition (a) locates the vertex of the source of the non-stationary points 
corresponding to the particular value of $d$.  The second condition identifies
non-stationary points that move from vertex $w$ to vertex $v$.

We can reduce condition (b) by first
simplifying the equality in (b) to see that $d = v(y) - t$.  
Because $0 < d$ we must have $0 < v(y) - t$, or equivalently, $t + 1 \le v(y)$.
With the inequality from (b) we have
\begin{equation}\label{cond b}
t + 1 \le v(y) \le a.
\end{equation}
Condition (a) translates to the inequality 
\begin{equation}\label{cond a}
s + d < 0 < t + 2d.
\end{equation}
These two inequalities have slightly different implications
depending on the type of $v$.

When $v$ is type $1^{0}$ we have $s < 0 < t$ so the right side of
(\ref{cond a}) is trivial since $0 < d$.  The left side
is equivalent to $s + v(y) - t < 0$ or $v(y) < t - s$.
Combining this with (\ref{cond b}) results in
\begin{equation}\label{cond type 1}
	t + 1 \le v(y) \le \min( a, t - s - 1).
\end{equation}

Otherwise, $v$ is type $2^{0}$ or $3^{0}$ and so
$s < t \le 0$.  The left side of (\ref{cond b}) is
trivial since $y \in \qtate{\mi}{\mi^{t}}$, so that
inequality can be replaced with $1 \le v(y) \le a$.
Replacing $d$ with $v(y) - t$ changes (\ref{cond a}) to
$s + v(y) - t < 0 < t + 2(v(y) - t)$.  The right side
is again trivial since $v(y) > 0$ and $t \le 0$.  To 
satisfy the left side, we need $v(y) < t - s$.  Combining
this with the modified inequality from (\ref{cond b})
gives
\begin{equation}\label{cond type 2,3}
	1 \le v(y) \le \min( a, t - s - 1).
\end{equation}

Provided $y$ satisfies inequality (\ref{cond type 1}) when
$v$ is type $1^{0}$ or inequality (\ref{cond type 2,3}) when
$v$ is type $2^{0}$ or $3^{0}$, there are
non-stationary points in the $I$-orbit of 
$w = (s + d, t + 2d)$ that move to the $I^{a}$-orbit
of $v$.  For all the $y$ of a fixed allowable valuation,
the non-stationary points in $S^{a}_{v}$ are enumerated by
\begin{eqnarray*}
	i & \in & \qtate{\ri}{\mi^{d-(s+d)}} = \qtate{\ri}{\mi^{-s}} \\
	j & \in & \qtate{\mi^{d-(t+2d)}}{\mi^{2d-(t+2d)}} = \qtate{\mi^{-v(y)}}{\mi^{-t}} \\
	z & \in & \qtate{\mi}{\mi^{(t+2d)-(s+d)-d}} = \qtate{\mi}{\mi^{t-s}}
\end{eqnarray*}
where $v(j) = -v(y)$ since $j = 1/y$.  Thus the non-stationary
points correspond to $-1 \ge v(j) \ge -\min( a, t - s - 1)$ when
$v$ is type $2^{0}$ or $3^{0}$ and $-(t+1) \ge v(j) \ge -\min( a, t - s - 1)$ when
$v$ is type $1^{0}$.  These points fit very nicely with the stationary
points at $v$ to give $S^{a}_{v}$ enumerated by
\begin{eqnarray*}
	i & \in & \qtate{\ri}{\mi^{-s}} \\
	j & \in & \qtate{\mi^{-x}}{\mi^{-t}} \\
	z & \in & \qtate{\mi}{\mi^{t-s}}
\end{eqnarray*}
where $x = \min(a,t-s-1)$ regardless of whether
$v$ is type $1^{0}$, $2^{0}$, or $3^{0}$.

The results from this section are summarized by
\begin{lemma}\label{sav enum}
Let $v=(s,t)$ be type $1^{0}$, $2^{0}$, or $3^{0}$ and let $x = \min(a,t-s-1)$.
The points in $S^{a}_{v}$ are enumerated by
\begin{center}
\begin{tabular}{|c|c|c|c|c|}\hline
\mycell{Type $v$} & \mycell{$i$} & \mycell{$j$} & \mycell{$y$} & \mycell{$z$}  \\ \hline
\mycell{ $3^{a}$} & \mycell{ $i \in \qtate{\ri}{\mi^{-s}}$} & \mycell{$j \in \qtate{\mi^{-x}}{\mi^{-t}}$} & \mycell{$y=0$} & \mycell{$z \in \qtate{\mi}{\mi^{t-s}}$} \\ \hline
\mycell{ $2^{a}$} & \mycell{$i \in \qtate{\ri}{\mi^{-s}}$} & \mycell{$j=0$} & \mycell{$y=0$} & \mycell{$z \in \qtate{\mi}{\mi^{t-s}}$} \\ \hline
\mycell{ $1^{a}$} & \mycell{$i \in \qtate{\ri}{\mi^{-s}}$} & \mycell{$j=0$} & \mycell{$y \in \qtate{\mi^{a+1}}{\mi^{t}}$} & \mycell{$z \in \qtate{\mi}{\mi^{t-s}}$} \\ \hline
\end{tabular}
\end{center}
\end{lemma}

\subsubsection{$v$ is Type $5^{0}$, $6^{0}$, or $7^{0}$}
This section and the next parallel \ref{alt retract I} and \ref{alt retract II}.
As before, we first prove 

\begin{lemma}\label{567 movement}
Let $M$ be a point in the $I$-orbit of the vertex $v = (s,t)$.
\begin{itemize}
\item[(a)] If $v$ is type $5^{0}$ or $6^{0}$ then $M$ is stationary
\item[(b)] If $v$ is type $7^{0}$ and $v( x_{M}) > a$ then $M$ is stationary
\item[(c)] If $v$ is type $7^{0}$ and $v( x_{M}) \le a$ then $M$ is in the $I^{a}$-orbit
	of $w = (s - 2d, t-d)$ where $d = s - v(x_{M})$
\end{itemize}
\end{lemma}

\emph{Part (a) \& (b)}: $M$ is stationary if $M \in I \cap I^{a}$.
From \ref{alt retract I} we know that 
\[ I = \mymatrix{ \ri^{\times} & \ri & \ri \\ 
	\mi^{a+1} & \ri^{\times} & \ri \\ \mi^{a+1} & \mi & \ri^{\times}}. \]
It is easy to see that any $M$ representing a point in the $I$-orbit
of a vertex of type $5^{0}$, $6^{0}$, or $7^{0}$ with
$v(x_M) > a$ is in $I \cap I^{a}$ and so
$M$ is stationary.

\emph{Part (c)}:\   Since $0 < v(x_{M}) < s$, we 
have the distance $d$ between $v$ and $w$ greater than zero.
To show that  
$M$ is in the $I^{a}$-orbit of $w$, it suffices to find $M'$ in $I^{a}$ such
that 
\[ M' w \in M v K.\]
As in section \ref{alt retract I}, we replace $w$ with the equivalent element $\pi^{d}w$
and consider
\[ M^{-1} M' \in \mymatrix{ \mi^{-d} &  \mi^{d-s} & \mi^{-t} \\ \mi^{s-d} &  \mi^{d} & \mi^{s-t} \\  
\mi^{t-d} &  \mi^{t-s+d} & \ri }.\]
To compute the matrix form of the left side,  we need to know 
the standard form for $M'$, which is determined by $w$'s type relative
to $I^{a}$.   

\begin{lemma}\label{w is type 5a}
Suppose that $v=(s,t)$ is type $7^0$ and 
$1 \le v(x) \le \min(a, s-1)$.  Define $d = s - v(x)$.  Then 
$ w = (s - 2d, t - d) $ is type $5^a$ and either type $5^0$, $6^0$, or $7^0$.
\end{lemma}

\emph{Proof}.\  The inequalities
\[ (t - d) < 0 \quad\textrm{and}\quad (t-d) < (s - 2d)\]
imply $w$ is type $5^0$, $6^0$, or $7^0$
and
\[ (t - d) < (s - 2d) < a\]
implies that $w$ is type $5^{a}$.

By assumption, $v(x) \le \min(a, s-1)$
so $v(x) < s$ and thus $d \ge 1$.  Since $v$ is 
type $7^0$ we know $t < 0$.  Therefore the inequality $(t - d) < 0$ is clear.

The inequality $(t - d) < (s - 2d)$ is equivalent to $t < v(x)$
which follows from $t < 0 < 1 \le v(x)$.
Finally, the inequality $(s - 2d) < a$ is equivalent to $2v(x) < a + s$
which follows from $v(x) \le a$ and $v(x) \le s - 1 < s$.~\eop

We now compute $M^{-1}M'$.

Because $w$ is type $5^{a}$,
\[ M' =  \mymatrix{1 & i' & j' \\ 0 & 1 & k' \\ 0 & 0 & 1} \]
and so
\begin{eqnarray*}
 M^{-1}M' & = & \mymatrix{ 1 & 0 & -j \\ -x & 1 & xj-k \\ 0 & 0 & 1 } \mymatrix{1 & i' & j' \\ 0 & 1 & k' \\ 0 & 0 & 1} \\
	& = & \mymatrix{ 1 & i' & (j'-j) \\  	-x & (1-xi') & [x(j-j') + (k'-k)] \\
	0 & 0 & 1}.
\end{eqnarray*}
Therefore, we need
\begin{equation}\label{m hat m condition II}
	\mymatrix{ 1 & i' & (j'-j) \\  	-x & (1-xi') & [x(j-j') + (k'-k)] \\
	0 & 0 & 1}
 \subset 
\mymatrix{ \mi^{-d} &  \mi^{d-s} & \mi^{-t} \\ \mi^{s-d} &  \mi^{d} & \mi^{s-t} \\  \mi^{t-d} &  \mi^{t-s+d} & \ri }
\end{equation}
for $M'$ to represent the point $M$ in the $I^{a}$-orbit of $w$.

Because $M$ is in the $I$-orbit of a type $7^{0}$ vertex and non-stationary, we have 
\[ x \in \qtate{\mi}{\mi^{s}}, \quad k \in \qtate{\ri}{\mi^{s-t}}, \textrm{\quad and \quad} 
 j \in \qtate{\ri}{\mi^{-t}} \]
with $v(x) \le \min(a, s-1)$.  To determine $M'$, first 
express $k$ as a polynomial in $\pi$
\[ k = k_{0} + k_{1} \pi^{1} + k_{2} \pi^{2} + \dots + k_{s-t-1} \pi^{s-t-1}. \]
and define
\[ \lceil k \rceil = k_{s-t-d} \pi^{s-t-d} + \dots + k_{s-t-1} \pi^{s-t-1}. \]
Let
\[\begin{array}{rclcr}  k' & = & k - \lceil k \rceil &\in& \qtate{\ri}{\mi^{s-t-d}}\\  
j' & = &\displaystyle j - \frac{\lceil k \rceil}{x}  &\in& \qtate{\ri}{\mi^{d-t}}\\  
i' & = &\displaystyle \frac{1}{x} &\in& \qtate{\mi^{d-s}}{\mi^{2d-s}}.
\end{array}\]

To verify that $ M^{-1}M' \in \pi^{-d}vKw^{-1}$, we check the 
interesting matrix entries in equation (\ref{m hat m condition II}). 
\begin{center}
\begin{tabular}{|c|c|}\hline
Entry & Check \\ \hline
$(1,2)$ & \mycell{$\displaystyle v(i') = -v(x) = d-s \textrm{ so } i' \in \mi^{d-s}$} \\ \hline
$(1,3)$ & \mycell{$\displaystyle v(j'-j) = v\left( \frac{\lceil k \rceil}{x} \right) = s- t-d + (d-s) = -t$} \\ \hline
$(2,1)$ & \mycell{$v(-x) = s - d$ so $-x \in \mi^{s-d}$} \\ \hline
$(2,2)$ & \mycell{$\begin{array}{lcl}\displaystyle  (1 - xi') 
    & = & \displaystyle  1 -  x\left( \frac{1}{x} \right) \\  
    & = &  0\end{array}$} \\ \hline
$(2,3)$ & \mycell{$ \begin{array}{lcl} \displaystyle x(j - j') + (k' - k) 
    & = & \displaystyle x\left( \frac{\lceil k \rceil}{x}\right) - \lceil k \rceil \\    
    & = & 0\end{array}$} \\ \hline
\end{tabular}
\end{center}

Thus the point $M$ in the $I$-orbit of $v$ is also represented by $M'$ in
the $I^{a}$-orbit of $w$.~\eop

We showed that for a non-stationary 
point the corresponding $w$ was of type $5^{0}$, $6^{0}$, or $7^{0}$.  
This shows that when we rearrange elements of $T$ by $I^{a}$-orbits, any 
non-stationary point is in the $I^{a}$-orbit of a vertex in $T$.

\begin{corollary}
\[ T = \bigsqcup_{\textrm{$v$ type $5^0$, $6^0$, or $7^0$}} T^{a}_{v} \]
\end{corollary}

\emph{Proof}.  It is clear from the definition of $T^{a}_{v}$
that $T^{a}_{v} \subset T$.

Let $M \in T$.  Then for some $v$ of type $5^0$, $6^0$, or $7^0$
\[ M \in \qtate{IvK}{K}.\]

If $M$ is stationary, then $M \in T^{a}_{v}$.  Otherwise, $M$ 
is non-stationary and there is a $w$ of type $5^0$, $6^0$, or $7^0$
such that $M \in T^{a}_{w}$.~\eop

\subsubsection{Structure of $T^{a}_{v}$ for $v$ Type $5^{0}$, $6^{0}$, or $7^{0}$}

In the previous section, we categorized the points of $T$ as stationary
and non-stationary.  Given a vertex $v$ of type $5^{0}$, $6^{0}$, or $7^{0}$,
we already know which points in $T^{a}_{v}$ are stationary.  To complete
the description of $T^{a}_{v}$, we need to determine the 
non-stationary points in the set.

When $v$ is type $7^{a}$ or $6^{a}$ the only points in $T^{a}_{v}$
are stationary points because non-stationary move to the $I^{a}$-orbit of 
type $5^{a}$ vertices.  
Since only type $7^{0}$ vertices can be type $7^{a}$ or $6^{a}$,
we can simply apply the stationary condition
$v(x) > a$ to the elements in the enumeration of points in the $I$-orbit of a
vertex of type $7^{0}$ to find $T^{a}_{v}$ in this case.  Therefore, 
when $v$ is type $7^{a}$ or $6^{a}$, we have 
$t < a \le s$ and $T^{a}_{v}$ is
enumerated by
\begin{eqnarray*}
	j & \in & \qtate{\ri}{\mi^{-t}} \\
	k & \in & \qtate{\ri}{\mi^{s-t}} \\
	x & \in & \qtate{\mi^{a+1}}{\mi^{s}}.
\end{eqnarray*}

Now we consider the case when $v$ is type $5^{a}$.
To list the stationary points, first suppose that $v$ is 
type $5^{a}$ and either type $7^{0}$ or type $6^{0}$.
Then $t < 0 \le s < a$.  If $v$ is type $6^0$, then $x = 0$.
Otherwise the stationary condition $v(x) > a$
applied to $x \in \qtate{\mi}{\mi^{s}}$ forces $x = 0$.  So 
the stationary points in $T^{a}_{v}$ are given by
\begin{eqnarray*}
	j & \in & \qtate{\ri}{\mi^{-t}} \\
	k & \in & \qtate{\ri}{\mi^{s-t}}. 
\end{eqnarray*}

When $v$ is type $5^{a}$ and $5^{0}$, 
$t < s < 0 < a$ and all the points are
stationary
\begin{eqnarray*}
	i & \in & \qtate{\ri}{\mi^{-s}} \\
	j & \in & \qtate{\ri}{\mi^{-t}} \\
	k & \in & \qtate{\ri}{\mi^{s-t}}. 
\end{eqnarray*}

For each positive integer $d$, there are non-stationary
points in $T$ that contribute to the $I^{a}$-orbit 
of $v = (s,t)$ if $v$ is type $5^{a}$
and
\begin{itemize}
\item[(a)] $w = (s + 2d, t + d)$ is a type $7^{0}$
\item[(b)] there exists $M \in \qtate{IwK}{K}$ so that $v(x) > a$ and $d = (s + 2d) - v(x)$ (where $x \equiv x_{M}$).
\end{itemize}
Condition (a) locates the vertex of the source of the non-stationary points 
corresponding to the particular value of $d$.  The second condition identifies
non-stationary points that move from vertex $w$ to vertex $v$.

We can reduce condition (b) by first
simplifying the equality in (b) to see that $d = v(x) - s$.  
Because $0 < d$ we must have $0 < v(x) - s$, or equivalently, $s + 1 \le v(x)$.
With the inequality from (b) we have
\begin{equation}\label{cond b II}
s + 1 \le v(x) \le a.
\end{equation}
Condition (a) translates to the inequality 
\begin{equation}\label{cond a II}
t + d < 0 < s + 2d.
\end{equation}
These two inequalities have slightly different implications
depending on the type of $v$.

When $v$ is type $7^{0}$ we have $t < 0 < s$ so the right side of
(\ref{cond a II}) is trivial since $0 < d$.  The left side
is equivalent to $t + v(x) - s < 0$ or $v(x) < s - t$.
Combining this with (\ref{cond b II}) results in
\begin{equation}\label{cond type 7}
	s + 1 \le v(x) \le \min( a, s - t - 1).
\end{equation}

Otherwise, $v$ is type $6^{0}$ or $5^{0}$ and so
$t < s \le 0$.  The left side of (\ref{cond b II}) is
trivial since $x \in \qtate{\mi}{\mi^{s}}$, so that
inequality can be replaced with $1 \le v(x) \le a$.
Replacing $d$ with $v(x)-s$ changes (\ref{cond a II}) to
$t + v(x) - s  < 0 < s + 2(v(x) - s)$.  The right side
is again trivial since $v(x) > 0$ and $s \le 0$.  To 
satisfy the left side, we need $v(x) < s - t$.  Combining
this with the modified inequality from (\ref{cond b II})
gives
\begin{equation}\label{cond type 6,5}
	1 \le v(x) \le \min( a, s - t - 1).
\end{equation}

Provided $x$ satisfies inequality (\ref{cond type 7}) when
$v$ is type $7^{0}$ or inequality (\ref{cond type 6,5}) when
$v$ is type $6^{0}$ or $5^{0}$, there are
non-stationary points in the $I$-orbit of 
$w = (s + 2d, t + d)$ that move to the $I^{a}$-orbit
of $v$.  For all the $x$ of a fixed allowable valuation,
the non-stationary points in $T^{a}_{v}$ are enumerated by
\begin{eqnarray*}
	k & \in & \qtate{\ri}{\mi^{(s+2d)-(t+d)-d}} = \qtate{\ri}{\mi^{s-t}} \\
	j & \in & \qtate{\ri}{\mi^{d-(t+d)}} = \qtate{\mi}{\mi^{-t}} \\
	i & \in & \qtate{\mi^{d-(s+2d)}}{\mi^{2d-(s+2d)}} = \qtate{\mi^{-v(x)}}{\mi^{-s}} 
\end{eqnarray*}
where $v(i) = -v(x)$ since $i = 1/x$.  Thus the non-stationary
points correspond to $-1 \ge v(i) \ge -\min( a, s - t - 1)$ when
$v$ is type $6^{0}$ or $5^{0}$ and $-(s+1) \ge v(i) \ge -\min( a, s - t - 1)$ when
$v$ is type $7^{0}$.  These points fit very nicely with the stationary
points at $v$ to give $T^{a}_{v}$ enumerated by
\begin{eqnarray*}
	k & \in & \qtate{\ri}{\mi^{s-t}} \\
	j & \in & \qtate{\ri}{\mi^{-t}} \\
	i & \in & \qtate{\mi^{-x}}{\mi^{-s}}
\end{eqnarray*}
where $x = \min(a,s-t-1)$ regardless of whether
$v$ is type $5^{0}$, $6^{0}$, or $7^{0}$.

The results from this section are summarized by
\begin{lemma}\label{tav enum}
Let $v=(s,t)$ be type $5^{0}$, $6^{0}$, or $7^{0}$ and let $x = \min(a,s-t-1)$.
The points in $T^{a}_{v}$ are enumerated by
\begin{center}
\begin{tabular}{|c|c|c|c|c|}\hline
\mycell{Type $v$} & \mycell{$i$} & \mycell{$j$} & \mycell{$k$} & \mycell{$x$}  \\ \hline
\mycell{ $5^{a}$} & 
	\mycell{ $i \in \qtate{\mi^{-x}}{\mi^{-s}}$} & 
	\mycell{ $j \in \qtate{\ri}{\mi^{-t}}$} & 
	\mycell{ $k \in \qtate{\ri}{\mi^{s-t}}$} & 
	\mycell{ $x=0$} \\ \hline
\mycell{ $6^{a}$} & 
	\mycell{ $i=0$} & 
	\mycell{ $j \in \qtate{\ri}{\mi^{-t}}$} & 
	\mycell{ $k \in \qtate{\ri}{\mi^{s-t}}$} & 
	\mycell{ $x=0$} \\ \hline
\mycell{ $7^{a}$} & 
	\mycell{ $i=0$} & 
	\mycell{ $j \in \qtate{\ri}{\mi^{-t}}$} & 
	\mycell{ $k \in \qtate{\ri}{\mi^{s-t}}$} & 
	\mycell{ $x \in \qtate{\mi^{a+1}}{\mi^{s}}$} \\ \hline
\end{tabular}
\end{center}
\end{lemma}

\section{Proof of the Main Theorem}\label{fixed point calculations}

In this section, we show that the sets 
$X^{\gamma} \cap S^{a}_{v}$, $X^{\gamma} \cap T^{a}_{v}$, and $X^{\gamma} \cap V^{0}_{v}$
are affine spaces by explicit calculation.
Since the standard form of the matrices representing points in $X$ 
vary depending on the type of $v$, we proceed case by case
through the twelve types.  The next section outlines the general approach and the following
sections applies the technique to specific regions.

\subsection{Fix Point Conditions}

Let $v=(s,t)$ and $M$ be a point in $\qtate{IvK}{K} = \qtate{I}{I_{(s,t)}}$.
$M$ is in $X^{\gamma}$ if $\gamma M = M$.
Since $\gamma \in I_{(s,t)}$, we can replace $\gamma M$
with $\gamma M \gamma^{-1}$ which has the advantage of
preserving the standard form of $M$.
Thus $M$ is fixed if and only if
$\gamma M \gamma^{-1} \in MI_{(s,t)}$ or equivalently if
\begin{equation}\label{fixed point}
M^{-1}\gamma M \gamma^{-1} \in  I_{(s,t)}.
\end{equation}
Immediately, we see that any vertex in the main apartment is fixed since $M$ is the identity in that case.

The matrix form of (\ref{fixed point}) depends upon the standard 
form of $M$ and thus $v$'s type.  For each of the six vertex types in $V$,
we write $M$ in standard form for that vertex type and compute the 
left side of (\ref{fixed point}).  We then analyze the resulting 
expression to determine $X^{\gamma} \cap V^{0}_{v}$.
In the cases when $v$ is in $S^{a}_{v}$ or $T^{a}_{v}$ we 
proceed in a similar fashion except the sets  
$I$ and $I_{(s,t)}$ are replaced with $I^{a}$ and $I^{a}_{(s,t)}$
because $S^{a}_{v}$ and $T^{a}_{v}$ are affine subsets of 
the $I^{a}$-orbit of $v$.

\subsection{$X^{\gamma} \cap V^{0}_{v}$}

The set $V$ contains vertices of type $4^{0}$,
$8^{0}$, $9^{0}$, $10^{0}$, $11^{0}$, and $12^{0}$.
In the following sections, the fixed point condition
is applied to the $I$-orbits of these vertex types 
and the dimension of the resulting affine space is determined.

\subsubsection{ $v$ type $4^{0}$}

The fixed point condition (\ref{fixed point}) is
\begin{eqnarray*}
	\mymatrix{ 1 & -i & -j \\ 0 & 1 & 0 \\ 0 & 0 & 1 } & \times &
	\mymatrix{ 1 & i\left( \frac{u_{1}}{u_{2}}\right) & j\left( \frac{u_{1}}{u_{3}}\right) \\
				0 & 1 & 0 \\ 0 & 0 & 1 } \\
	 & = & \mymatrix{ 1 & i\left( \frac{u_{1}}{u_{2}} - 1\right) & j\left( \frac{u_{1}}{u_{3}} - 1\right) \\
				0 & 1 & 0 \\ 0 & 0 & 1 } \\
 	& \in & \mymatrix{ 1 & \mi^{-s} & \mi^{-t} \\ 0 & 1 & 0 \\ 0 & 0 & 1 }.
\end{eqnarray*}
Thus $M$ is fixed if and only if $v(i) \ge -s - m$ and $v(j) \ge -t - m$.  

Since $i \in \qtate{\ri}{\mi^{-s}}$ and $j \in \qtate{\ri}{\mi^{-t}}$ 
when $v$ is type $4^{0}$, the set $X^{\gamma} \cap V^{0}_{v}$ is an affine 
space of dimension $\min(m, -s) + \min(m, -t)$.  The variable $i$ 
contributes $\min(m, -s)$ to the dimension and the variable $j$ 
contributes $\min(m, -t)$ to the dimension.

\subsubsection{ $v$ type $8^{0}$}

The fixed point condition (\ref{fixed point}) is
\begin{eqnarray*}
	\mymatrix{ 1 & 0 & 0 \\ -x & 1 & -k \\ 0 & 0 & 1 } & \times &
	\mymatrix{ 1 & 0 & 0 \\
				x\left( \frac{u_{2}}{u_{1}}\right) & 1 & k\left( \frac{u_{2}}{u_{3}}\right) \\ 0 & 0 & 1 } \\
	 & = & \mymatrix{ 1 & 0 & 0 \\
				x\left( \frac{u_{2}}{u_{1}} - 1\right) & 1 & k\left( \frac{u_{2}}{u_{3}} - 1\right) \\ 0 & 0 & 1 } \\
 	& \in & \mymatrix{ 1 & 0 & 0 \\ \mi^{s} & 1 & \mi^{s-t} \\ 0 & 0 & 1 }.
\end{eqnarray*}
Thus $M$ is fixed if and only if $v(x) \ge s - m$ and $v(k) \ge s-t - n$.  

Since $x \in \qtate{\mi}{\mi^{s}}$ and $k \in \qtate{\ri}{\mi^{s-t}}$
when $v$ is type $8^{0}$, the set $X^{\gamma} \cap V^{0}_{v}$ is an affine 
space of dimension $\min(m, s-1) + \min(n, s-t)$.

\subsubsection{ $v$ type $9^{0}$}\label{v type 90}
 
The fixed point condition (\ref{fixed point}) is
\begin{eqnarray*}
	\mymatrix{ 1 & 0 & 0 \\ ky - x & 1 & -k \\ -y & 0 & 1 } & \times &
	\mymatrix{ 1 & 0 & 0 \\ 
					x\left( \frac{u_{2}}{u_{1}} \right) & 1 & k\left( \frac{u_{2}}{u_{3}}\right) \\
					y\left( \frac{u_{3}}{u_{1}}\right) & 0 & 1} \\
	 & = & \mymatrix{ 1 & 0 & 0 \\
				 ky\left( 1-\frac{u_{3}}{u_{1}}\right) + x\left( \frac{u_{2}}{u_{1}} - 1\right)
					& 1 & k\left( \frac{u_{2}}{u_{3}} - 1\right) \\
				 y\left( \frac{u_{3}}{u_{1}} - 1\right) & 0 & 1} \\
 	& \in & \mymatrix{ 1 & 0 & 0 \\ \mi^{s} & 1 & \mi^{s-t} \\ \mi^{t} & 0 & 1 }.
\end{eqnarray*}
From the enumeration of the points in the $I$-orbit of $v$, we 
have $k \in \qtate{\ri}{\mi^{s-t}}$, $x\in \qtate{\mi}{\mi^{s}}$, and $y \in \qtate{\mi}{\mi^{t}}$.
Incorporating the fixed point conditions, we then have
\begin{eqnarray*}
s-t > &v(k)& \ge \max(s-t-n,0) \\
t > &v(y)& \ge \max(t-m,1)
\end{eqnarray*}
and
\[ x\in \qtate{\mi}{\mi^{s}} \textrm{ such that } ky\left( 1-\frac{u_{3}}{u_{1}}\right) + x\left( \frac{u_{2}}{u_{1}} - 1\right) \in \mi^{s}. \]

Given a $k$ and $y$ that satisfy the above inequalities, 
$x$ is partially determined by those values.  We write 
$x = x' + x''$ where $x'$ is the determined portion
and $x''$ is the free portion.
Since \[ v\left( ky\left( 1-\frac{u_{3}}{u_{1}}\right) \right) > 0 + 1 + m \]
and
\[  v\left( x\left( \frac{u_{2}}{u_{1}} - 1\right) \right) \ge 1 + m  \]
we can always set 
\[ x' = -ky\frac{ \left( 1-\frac{u_{3}}{u_{1}}\right) }{ \left( \frac{u_{2}}{u_{1}} - 1\right) }
 \in \qtate{\mi}{ \mi^{s}} \]
and then
\[ ky\left( 1-\frac{u_{3}}{u_{1}}\right) + x'\left( \frac{u_{2}}{u_{1}} - 1\right) \in \mi^{s}. \]
For any $x'' \in \qtate{ \mi}{\mi^{s}}$ such that
\[ v(x'') \ge \max( s - m, 1) \]
we have 
\[ x'' \left( \frac{u_{2}}{u_{1}} - 1\right) \in \mi^{s} \]
so for $x = x' + x''$, 
\[ ky\left( 1-\frac{u_{3}}{u_{1}}\right) + x\left( \frac{u_{2}}{u_{1}} - 1\right) \in \mi^{s} \]
which satisfies the fixed point condition.

Therefore, when $v$ is type $9^{0}$, the set $X^{\gamma} \cap V^{0}_{v}$ is an affine 
space of dimension $\min( n, s-t) + \min( m, t-1) + \min( m, s-1)$.

\subsubsection{ $v$ type $10^{0}$}

The fixed point condition (\ref{fixed point}) is
\begin{eqnarray*}
	\mymatrix{ 1 & 0 & 0 \\ -x & 1 & 0 \\ -y & 0 & 1 } & \times &
	\mymatrix{ 1 & 0 & 0 \\
				x\left( \frac{u_{2}}{u_{1}}\right) & 1 & 0 \\ y\left( \frac{u_{3}}{u_{1}}\right) & 0 & 1 } \\
	 & = & \mymatrix{ 1 & 0 & 0 \\
				x\left( \frac{u_{2}}{u_{1}} - 1\right) & 1 & 0 \\ y\left( \frac{u_{3}}{u_{1}} - 1\right) & 0 & 1 } \\
 	& \in & \mymatrix{ 1 & 0 & 0 \\ \mi^{s} & 1 & 0 \\ \mi^{t} & 0 & 1 }.
\end{eqnarray*}
Thus $M$ is fixed if and only if $v(x) \ge s - m$ and $v(y) \ge t - m$.  

Since $x \in \qtate{\mi}{\mi^{s}}$ and $y \in \qtate{\mi}{\mi^{t}}$
when $v$ is type $10^{0}$, the set $X^{\gamma} \cap V^{0}_{v}$ is an affine 
space of dimension $\min(m, s-1) + \min(m, t-1)$.

\subsubsection{ $v$ type $11^{0}$}

The fixed point condition (\ref{fixed point}) is
\begin{eqnarray*}
	\mymatrix{ 1 & 0 & 0 \\ -x & 1 & 0 \\ zx - y & -z & 1 } & \times & 
	\mymatrix{ 1 & 0 & 0 \\ 
					x\left( \frac{u_{2}}{u_{1}} \right) & 1 & 0 \\
					y\left( \frac{u_{3}}{u_{1}}\right) & z\left( \frac{u_{3}}{u_{2}}\right) & 1} \\
	 & = & \mymatrix{ 1 & 0 & 0 \\
				x\left( \frac{u_{2}}{u_{1}} - 1\right)
					& 1 & 0 \\
				  zx\left( 1-\frac{u_{2}}{u_{1}}\right) + y\left( \frac{u_{3}}{u_{1}} - 1\right) & 
				z\left( \frac{u_{3}}{u_{2}} - 1\right) & 1} \\
 	& \in & \mymatrix{ 1 & 0 & 0 \\ \mi^{s} & 1 & 0 \\ \mi^{t} & \mi^{t-s} & 1 }.
\end{eqnarray*}
From the enumeration of the points in the $I$-orbit of $v$, we 
have $z \in \qtate{\mi}{\mi^{t-s}}$, $x\in \qtate{\mi}{\mi^{s}}$, and $y \in \qtate{\mi}{\mi^{t}}$.
Incorporating the fixed point conditions, we then have
\begin{eqnarray*}
t-s > &v(z)& \ge \max(t-s-n,1) \\
s > &v(x)& \ge \max(s-m,1)
\end{eqnarray*}
and
\[ y\in \qtate{\mi}{\mi^{t}} \textrm{ such that } zx\left( 1-\frac{u_{2}}{u_{1}}\right) + y\left( \frac{u_{3}}{u_{1}} - 1\right) \in \mi^{t}. \]

Just as $x$ in section \ref{v type 90} was considered part determined and part free, the entry $y$ 
is partially determined by the choice of $z$ and $x$.  First, suppose that we choose  suitable
$z$ and $x$ that satisfy the fix point condition.  Since 
\[ v\left( zx\left( 1-\frac{u_{2}}{u_{1}}\right) \right) \ge 1 + 1 + m = 2 + m\]
and
\[ v\left( y\left( \frac{u_{3}}{u_{1}} - 1\right) \right) \ge 1 + m \]
we can choose a suitable $y$ (in the same manner as we chose $x$ in section \ref{v type 90})
so that $zx\left( 1-\frac{u_{2}}{u_{1}}\right) + y\left( \frac{u_{3}}{u_{1}} - 1\right) \in \mi^{t}$
with $y''$, the free part of $y$, restricted by
\[ t > v(y'') \ge \max(t-m, 1).\]

Therefore, when $v$ is type $11^{0}$, the set $X^{\gamma} \cap V^{0}_{v}$ is an affine 
space of dimension $\min( n, t-s-1) + \min( m, s-1) + \min( m, t-1)$.

\subsubsection{ $v$ type $12^{0}$}

The fixed point condition (\ref{fixed point}) is
\begin{eqnarray*}
	\mymatrix{ 1 & 0 & 0 \\ 0 & 1 & 0 \\ -y & -z & 1 } & \times &
	\mymatrix{ 1 & 0 & 0 \\
				0 & 1 & 0 \\ y\left( \frac{u_{3}}{u_{1}}\right) & z\left( \frac{u_{3}}{u_{2}}\right) & 1 } \\
	 & = & \mymatrix{ 1 & 0 & 0 \\
				0 & 1 & 0 \\ y\left( \frac{u_{3}}{u_{1}} - 1\right) & z\left( \frac{u_{3}}{u_{2}} - 1\right) & 1 } \\
 	& \in & \mymatrix{ 1 & 0 & 0 \\ 0 & 1 & 0 \\ \mi^{t} & \mi^{t-s} & 1 }.
\end{eqnarray*}
Thus $M$ is fixed if and only if $v(z) \ge t- s - n$ and $v(y) \ge t - m$.  
Since $z \in \qtate{\mi}{\mi^{t-s}}$ and $y \in \qtate{\mi}{\mi^{t}}$
when $v$ is type $12^{0}$, the set $X^{\gamma} \cap V^{0}_{v}$ is an affine 
space of dimension $\min(n, t-s-1) + \min(m, t-1)$.

\subsection{$X^{\gamma} \cap S^{a}_{v}$}

The set $S$ contains vertices of type $1^{a}$,
$2^{a}$, and $3^{a}$.
In the following sections, the fixed point conditions
are applied to the $I^{a}$-orbits of these vertex types 
and the dimension of the resulting affine space is determined.

Throughout, we reference the enumeration of points
in $S^{a}_{v}$ described in Lemma~\ref{sav enum}.

\subsubsection{ $v$ is type $1^{a}$}

The fixed point condition (\ref{fixed point}) is
\begin{eqnarray*}
	\mymatrix{ 1 & -i & 0 \\ 0 & 1 & 0 \\ -y & yi-z & 1 } & \times &
	\mymatrix{ 1 & i\left( \frac{u_{1}}{u_{2}} \right) & 0 \\ 
					0 & 1 & 0 \\
					y\left( \frac{u_{3}}{u_{1}}\right) & z\left( \frac{u_{3}}{u_{2}}\right) & 1} \\
	 & = & \mymatrix{ 1 & i\left( \frac{u_{1}}{u_{2}} - 1\right) & 0 \\
					 0 & 1 & 0 \\
					 y\left( \frac{u_{3}}{u_{1}} - 1\right) & 
						yi\left( 1 - \frac{u_{1}}{u_{2}}\right) + z\left( \frac{u_{3}}{u_{2}} - 1\right)
					 & 1} \\
 	& \in & \mymatrix{ 1 & \mi^{-s} & 0 \\ 0 & 1 & 0 \\ \mi^{t} & \mi^{t-s} & 1 }.
\end{eqnarray*}

Since $i \in \qtate{\ri}{\mi^{-s}}$ and $y \in \qtate{\mi^{a+1}}{\mi^{t}}$, to satisfy the 
fixed point condition we require that
\begin{eqnarray*}
-s > & v(i) & \ge \max( -s-m, 0) \\
t > & v(y) & \ge \max(t-m, a+1)
\end{eqnarray*}
and 
\[ z \in \qtate{\mi}{\mi^{t-s}} \textrm{ such that } 
	yi\left( 1 - \frac{u_{1}}{u_{2}}\right) + z\left( \frac{u_{3}}{u_{2}} - 1\right) \in \mi^{t-s}.\]

Just as $x$ in section \ref{v type 90} was considered part determined and part free, the entry $z$ 
is partially determined by the choice of $i$ and $y$.  First, suppose that we choose  suitable
$i$ and $y$ that satisfy the fix point condition.  Since 
\[ v\left( yi \left( 1 - \frac{u_{1}}{u_{2}} \right) \right) \ge (a+1) + 0 + m = n + 1\]
and
\[ v\left( z\left( \frac{u_{3}}{u_{2}} - 1\right) \right) \ge 1 + n \]
we can choose a suitable $z$ (in the same manner as we chose $x$ in section \ref{v type 90})
so that $yi\left( 1 - \frac{u_{1}}{u_{2}}\right) + z\left( \frac{u_{3}}{u_{2}} - 1\right) \in \mi^{t-s}$
with $z''$, the free part of $z$, restricted by
\[ t - s > v(z'') \ge \max(t-s-n, 1).\]

Therefore, when $v$ is in $S$ and type $1^{a}$ the set $X^{\gamma} \cap S^{a}_{v}$ 
is an affine space of dimension $\min(m, -s) + \min(m, t-(a+1)) + \min(n,t-s-1)$.

\subsubsection{ $v$ is type $2^{a}$}

The fixed point condition (\ref{fixed point}) is
\begin{eqnarray*}
	\mymatrix{ 1 & -i & 0 \\ 0 & 1 & 0 \\ 0 & -z & 1 } & \times &
	\mymatrix{ 1 & i\left( \frac{u_{1}}{u_{2}}\right) & 0 \\
				0 & 1 & 0 \\ 0 & z\left( \frac{u_{3}}{u_{2}}\right) & 1 } \\
	 & = & \mymatrix{ 1 & i\left( \frac{u_{1}}{u_{2}} - 1\right) & 0 \\
				0 & 1 & 0 \\ 0 & z\left( \frac{u_{3}}{u_{2}} - 1\right) & 1 } \\
 	& \in & \mymatrix{ 1 & \mi^{-s} & 0 \\ 0 & 1 & 0 \\ 0 & \mi^{t-s} & 1 }.
\end{eqnarray*}
Thus $M$ is fixed if and only if $v(i) \ge -s - m$ and $v(z) \ge t-s - n$.  

Since $i \in \qtate{\ri}{\mi^{-s}}$ and $z \in \qtate{\mi}{\mi^{t-s}}$
when $v$ is type $2^{a}$, the set $X^{\gamma} \cap S^{a}_{v}$ is an affine 
space of dimension $\min(m, -s) + \min(n, t-s-1)$.

\subsubsection{ $v$ is type $3^{a}$}

The fixed point condition (\ref{fixed point}) is
\begin{eqnarray*}
	\mymatrix{ 1 & jz-i & -j \\ 0 & 1 & 0 \\ 0 & -z & 1 } & \times &
	\mymatrix{ 1 & i\left( \frac{u_{1}}{u_{2}} \right) & j\left( \frac{u_{1}}{u_{3}} \right) \\ 
					0 & 1 & 0 \\
					0 & z\left( \frac{u_{3}}{u_{2}}\right) & 1} \\
	 & = & \mymatrix{ 1 & jz\left( 1 - \frac{u_{3}}{u_{2}}\right) + i\left( \frac{u_{1}}{u_{2}} - 1\right)
						  & j\left( \frac{u_{1}}{u_{3}} - 1\right) \\
					 0 & 1 & 0 \\
					 0 & z\left( \frac{u_{3}}{u_{2}} - 1\right) & 1} \\
 	& \in & \mymatrix{ 1 & \mi^{-s} & \mi^{-t} \\ 0 & 1 & 0 \\ 0 & \mi^{t-s} & 1 }.
\end{eqnarray*}

Since $j \in \qtate{\mi^{-x}}{\mi^{-t}}$ (where $x = \min(a,t-s-1)$) and 
$z \in \qtate{\mi}{\mi^{t-s}}$, to satisfy the 
fixed point condition we require that
\begin{eqnarray*}
-t > & v(j) & \ge \max( -t-m, -x) \\
t-s > & v(z) & \ge \max(t-s-n, 1)
\end{eqnarray*}
and 
\[ i \in \qtate{\ri}{\mi^{-s}} \textrm{ such that } 
	jz\left( 1 - \frac{u_{3}}{u_{2}}\right) + i\left( \frac{u_{1}}{u_{2}} - 1\right) \in \mi^{-s}.\]

Just as $x$ in section \ref{v type 90} was considered part determined and part free, the entry $i$ 
is partially determined by the choice of $j$ and $z$.  First, suppose that we choose a suitable
$j$ and $z$ that satisfies the fix point condition.  Since 
\begin{eqnarray*}
 v\left( jz\left( 1 - \frac{u_{3}}{u_{2}} \right) \right) & \ge & -\min(a,t-s-1) + 1 + n \\
& = & \max( -a, -(t-s-1)) + 1 + n \\
& \ge & -(n-m) + 1 + n \\
& = & m + 1
\end{eqnarray*}
and
\[ v\left( i\left( \frac{u_{1}}{u_{2}} - 1\right) \right) \ge m \]
we can choose a suitable $i$ (in the same manner as we chose $x$ in section \ref{v type 90})
so that $jz\left( 1 - \frac{u_{3}}{u_{2}}\right) + i\left( \frac{u_{1}}{u_{2}} - 1\right) \in \mi^{-s}$
with $i''$, the free part of $i$, restricted by
\[ -s > v(i'') \ge \max(-s-m, 0).\]

Therefore, when $v$ is in $S$ and type $3^{a}$ the set $X^{\gamma} \cap S^{a}_{v}$ 
is an affine space of dimension $\min(m, -s) + \min(n, t-s-1)) + \min(m,a-t,-s-1)$.

\subsection{$X^{\gamma} \cap T$}

The set $T$ contains vertices of type $5^{a}$,
$6^{a}$, and $7^{a}$.
In the following sections, the fixed point conditions
are applied to the $I^{a}$-orbits of these vertex types 
and the dimension of the resulting affine space is determined.

Throughout, we reference the enumeration of points
in $T^{a}_{v}$ described in Lemma~\ref{tav enum}.

\subsubsection{ $v$ is type $7^{a}$}

The fixed point condition (\ref{fixed point}) is
\begin{eqnarray*}
	\mymatrix{ 1 & 0 & -j \\ -x & 1 & xj-k \\ 0 & 0 & 1 } & \times &
	\mymatrix{ 1 & 0 & j\left( \frac{u_{1}}{u_{3}}\right) \\ 
					x\left( \frac{u_{2}}{u_{1}}\right) & 1 & k\left( \frac{u_{2}}{u_{3}}\right) \\
					0 & 0 & 1} \\
	 & = & \mymatrix{ 1 & 0 & j\left( \frac{u_{1}}{u_{3}} - 1\right) \\
					 x\left( \frac{u_{2}}{u_{1}} - 1\right) & 1 & 
						xj\left( 1 - \frac{u_{1}}{u_{3}}\right) + k\left( \frac{u_{2}}{u_{3}} - 1\right) \\
					 0 & 0 & 1} \\
 	& \in & \mymatrix{ 1 & 0 & \mi^{-t} \\ \mi^{s} & 1 & \mi^{s-t} \\ 0 & 0 & 1 }.
\end{eqnarray*}

Since $j \in \qtate{\ri}{\mi^{-t}}$ and $x \in \qtate{\mi^{a+1}}{\mi^{s}}$, to satisfy the 
fixed point condition we require that
\begin{eqnarray*}
-t > & v(j) & \ge \max( -t-m, 0) \\
s > & v(x) & \ge \max(s-m, a+1)
\end{eqnarray*}
and 
\[ k \in \qtate{\ri}{\mi^{s-t}} \textrm{ such that } 
	xj\left( 1 - \frac{u_{1}}{u_{3}}\right) + k\left( \frac{u_{2}}{u_{3}} - 1\right) \in \mi^{s-t}.\]

Just as $x$ in section \ref{v type 90} was considered part determined and part free, the entry $k$ 
is partially determined by the choice of $j$ and $x$.  First, suppose that we choose a suitable
$j$ and $x$ that satisfies the fix point condition.  Since 
\[ v\left( xj\left( 1 - \frac{u_{1}}{u_{3}}\right) \right) \ge (a+1) + 0 + m = n + 1\]
and
\[ v\left( k\left( \frac{u_{2}}{u_{3}} - 1\right) \right) \ge 1 + n \]
we can choose a suitable $k$ (in the same manner as we chose $x$ in section \ref{v type 90})
so that $xj\left( 1 - \frac{u_{1}}{u_{3}}\right) + k\left( \frac{u_{2}}{u_{3}} - 1\right) \in \mi^{s-t}$
with $k''$, the free part of $k$, restricted by
\[ s-t > v(k'') \ge \max(s-t-n, 0).\]

Therefore, when $v$ is in $T$ and type $7^{a}$ the set $X^{\gamma} \cap T^{a}_{v}$ 
is an affine space of dimension $\min(m, -t) + \min(m, s-(a+1)) + \min(n,s-t)$.

\subsubsection{ $v$ is type $6^{a}$}

\begin{eqnarray*}
	\mymatrix{ 1 & 0 & -j \\ 0 & 1 & -k \\ 0 & 0 & 1 }
	\mymatrix{ 1 & 0 & j\left( \frac{u_{1}}{u_{3}}\right) \\ 
					0 & 1 & k\left( \frac{u_{2}}{u_{3}}\right) \\
					0 & 0 & 1}
	 & = & \mymatrix{ 1 & 0 & j\left( \frac{u_{1}}{u_{3}} - 1\right) \\
					 0 & 1 & k\left( \frac{u_{2}}{u_{3}} - 1\right) \\
					 0 & 0 & 1} \\
 	& \in & \mymatrix{ 1 & 0 & \mi^{-t} \\ 0 & 1 & \mi^{s-t} \\ 0 & 0 & 1 }.
\end{eqnarray*}
Thus $M$ is fixed if and only if $v(j) \ge -t - m$ and $v(k) \ge s-t - n$.  

Since $j \in \qtate{\ri}{\mi^{-t}}$ and $k \in \qtate{\ri}{\mi^{s-t}}$
when $v$ is type $6^{a}$, the set $X^{\gamma} \cap T^{a}_{v}$ is an affine 
space of dimension $\min(m, -t) + \min(n, s-t)$.

\subsubsection{ $v$ is type $5^{a}$}

The fixed point condition (\ref{fixed point}) is
\begin{eqnarray*}
	\mymatrix{ 1 & -i & ki-j \\ 0 & 1 & -k \\ 0 & 0 & 1 } & \times &
	\mymatrix{ 1 & i\left( \frac{u_{1}}{u_{2}} \right) & j\left( \frac{u_{1}}{u_{3}} \right) \\ 
					0 & 1 & k\left( \frac{u_{2}}{u_{3}}\right) \\
					0 & 0 & 1} \\
	 & = & \mymatrix{ 1 & i\left( \frac{u_{1}}{u_{2}} - 1\right)
						  & ki\left( 1 - \frac{u_{2}}{u_{3}}\right) + j\left( \frac{u_{1}}{u_{3}} - 1\right) \\
					 0 & 1 & k\left( \frac{u_{2}}{u_{3}} - 1\right) \\
					 0 & 0 & 1} \\
 	& \in & \mymatrix{ 1 & \mi^{-s} & \mi^{-t} \\ 0 & 1 & \mi^{s-t} \\ 0 & 0 & 1 }.
\end{eqnarray*}

Since $i \in \qtate{\mi^{-x}}{\mi^{-s}}$ (where $x = \min(a,s-t-1)$) and 
$k \in \qtate{\ri}{\mi^{s-t}}$, to satisfy the 
fixed point condition we require that
\begin{eqnarray*}
-s > & v(i) & \ge \max( -s-m, -x) \\
s-t > & v(k) & \ge \max(s-t-n, 0)
\end{eqnarray*}
and 
\[ j \in \qtate{\ri}{\mi^{-t}} \textrm{ such that } 
	ki\left( 1 - \frac{u_{2}}{u_{3}}\right) + j\left( \frac{u_{1}}{u_{3}} - 1\right) \in \mi^{-t}.\]

Just as $x$ in section \ref{v type 90} was considered part determined and part free, the entry $j$ 
is partially determined by the choice of $k$ and $i$.  First, suppose that we choose a suitable
$k$ and $i$ that satisfies the fix point condition.  Since 
\begin{eqnarray*}
 v\left( ki\left( 1 - \frac{u_{2}}{u_{3}}\right) \right) & \ge & 0 -\min(a,s-t-1) + n \\
& = & \max( -a, -(s-t-1)) + n \\
& \ge & -(n-m) +  n \\
& = & m
\end{eqnarray*}
and
\[ v\left( j\left( \frac{u_{1}}{u_{3}} - 1\right) \right) \ge m \]
we can choose a suitable $j$ (in the same manner as we chose $x$ in section \ref{v type 90})
so that $ki\left( 1 - \frac{u_{2}}{u_{3}}\right) + j\left( \frac{u_{1}}{u_{3}} - 1\right) \in \mi^{-t}$
with $j''$, the free part of $j$, restricted by
\[ -t > v(j'') \ge \max(-t-m, 0).\]

Therefore, when $v$ is in $T$ and type $5^{a}$ the set $X^{\gamma} \cap T^{a}_{v}$ 
is an affine space of dimension $\min(m, -t) + \min(n, s-t)) + \min(m,a-s,-t-1)$.

\section{Topological Considerations}

In the previous sections, we developed a decomposition of 
$X^{\gamma}$ into affine pieces.  We now examine the topology 
of the decomposition to verify it is a paving by affine 
spaces.

\subsection{Locally Closed}

To show that the affine pieces of $X^{\gamma}$ are 
locally closed, it is enough to show that the sets
$V^{0}_{v}$, $S^{a}_{v}$, and $T^{a}_{v}$ are locally closed.
Every $V^{0}_{v}$ set is simply an $I$-orbit of a vertex and thus locally closed.

From Lemma~\ref{sav enum} we know that 
$S^{a}_{v}$ is a closed affine subspace of $\qtate{I^{a}vK}{K}$,
which is locally closed because it is the $I^{a}$-orbit of $v$,
and thus $S^{a}_{v}$ is locally closed.  Similarly, from Lemma~\ref{tav enum}
we can conclude that $T^{a}_{v}$ is locally closed.

\subsection{A Filtration by Closed Sets}

In this section, we exhibit a filtration of $X$ by closed subsets
that are formed by an increasing union of the affine pieces of 
$X$ described in the previous sections.  Since the affine pieces
are in a one-to-one correspondence with the vertices in the main apartment
of $X$, we can describe the filtration by placing an order on those
vertices.  If $v_{i} = (s,t)$ is the $i$-th vertex, let $\mathbb{A}_{i}$ denote 
the affine space containing $v_{i}$ and then
\[ \mathbb{A}_{0} \subset \mathbb{A}_{0} \cup \mathbb{A}_{1} \subset \dots \]
will give the filtration.

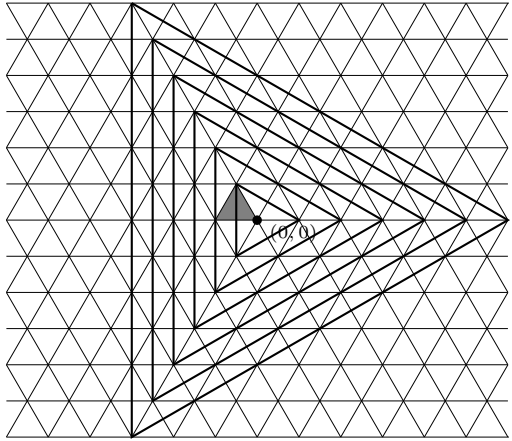
\begin{figure}
	\begin{center}
    \psset{xunit=0.10937499999999975 in}
    \psset{yunit=0.18944305707784575 in}
    \begin{pspicture}(0,0)(24,12)
    \scriptsize
	
	\pspolygon[fillstyle=solid,fillcolor=gray,linestyle=none](12,6)(11,7)(10,6)

	\psline(11,7)(11,5)(14,6)(11,7)	
	\psline(10,8)(10,4)(16,6)(10,8)	
	\psline(9,9)(9,3)(18,6)(9,9)
	\psline(8,10)(8,2)(20,6)(8,10)	
	\psline(7,11)(7,1)(22,6)(7,11)	
	\psline(6,12)(6,0)(24,6)(6,12)

        \psline[linewidth=0.1pt](0,2)(2,0)\psline[linewidth=0.1pt](0,4)(4,0)\psline[linewidth=0.1pt](0,6)(6,0)\psline[linewidth=0.1pt](0,8)(8,0)\psline[linewidth=0.1pt](0,10)(10,0)\psline[linewidth=0.1pt](0,12)(12,0)\psline[linewidth=0.1pt](2,12)(14,0)\psline[linewidth=0.1pt](4,12)(16,0)\psline[linewidth=0.1pt](6,12)(18,0)\psline[linewidth=0.1pt](8,12)(20,0)\psline[linewidth=0.1pt](10,12)(22,0)\psline[linewidth=0.1pt](12,12)(24,0)\psline[linewidth=0.1pt](14,12)(24,2)\psline[linewidth=0.1pt](16,12)(24,4)\psline[linewidth=0.1pt](18,12)(24,6)\psline[linewidth=0.1pt](20,12)(24,8)\psline[linewidth=0.1pt](22,12)(24,10)\psline[linewidth=0.1pt](24,12)(24,12)\psline[linewidth=0.1pt](0,10)(2,12)\psline[linewidth=0.1pt](0,8)(4,12)\psline[linewidth=0.1pt](0,6)(6,12)\psline[linewidth=0.1pt](0,4)(8,12)\psline[linewidth=0.1pt](0,2)(10,12)\psline[linewidth=0.1pt](0,0)(12,12)\psline[linewidth=0.1pt](2,0)(14,12)\psline[linewidth=0.1pt](4,0)(16,12)\psline[linewidth=0.1pt](6,0)(18,12)\psline[linewidth=0.1pt](8,0)(20,12)\psline[linewidth=0.1pt](10,0)(22,12)\psline[linewidth=0.1pt](12,0)(24,12)\psline[linewidth=0.1pt](14,0)(24,10)\psline[linewidth=0.1pt](16,0)(24,8)\psline[linewidth=0.1pt](18,0)(24,6)\psline[linewidth=0.1pt](20,0)(24,4)\psline[linewidth=0.1pt](22,0)(24,2)\psline[linewidth=0.1pt](24,0)(24,0)\psline[linewidth=0.1pt](0,0)(24,0)\psline[linewidth=0.1pt](0,1)(24,1)\psline[linewidth=0.1pt](0,2)(24,2)\psline[linewidth=0.1pt](0,3)(24,3)\psline[linewidth=0.1pt](0,4)(24,4)\psline[linewidth=0.1pt](0,5)(24,5)\psline[linewidth=0.1pt](0,6)(24,6)\psline[linewidth=0.1pt](0,7)(24,7)\psline[linewidth=0.1pt](0,8)(24,8)\psline[linewidth=0.1pt](0,9)(24,9)\psline[linewidth=0.1pt](0,10)(24,10)\psline[linewidth=0.1pt](0,11)(24,11)\psline[linewidth=0.1pt](0,12)(24,12)

\psdots[dotstyle=*](12,6) 
\uput[330](12,6){$(0,0)$}
	\end{pspicture}
    \end{center}
\caption{Triangular grouping of vertices}\label{triangles}
\end{figure}

The base point $v_{0} = (0,0)$ is a single point and the first closed subset
in our filtration.  To describe the order on the remaining vertices, we first coarsely group
the vertices into increasingly larger sets defined by triangular bounds as in figure \ref{triangles}.
The smallest triangle, $\Delta_{1}$, represents a set of four vertices
(three on the boundary and the base point $v_{0}$).  Triangle $\Delta_{2}$
is the next largest and it represents the set of six vertices on its
boundary plus the four vertices from $\Delta_{1}$.  We define $\Delta_{0} = \{ v_{0} \}$.

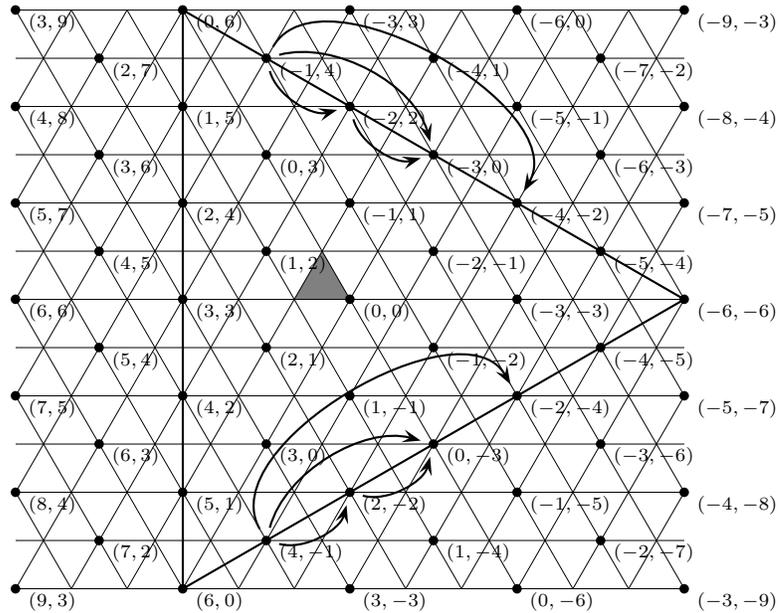
\begin{figure}[!hbt]
	\begin{center}
   \psset{xunit=0.145833333333333 in}
    \psset{yunit=0.252590742770461 in}
    \begin{pspicture}(0,0)(24,12)
    \scriptsize
\pspolygon[fillstyle=solid,fillcolor=gray,linestyle=none](12,6)(11,7)(10,6)
\psline(6,0)(6,12)(24,6)(6,0)

\pnode(9,11){A}
\pnode(12,10){B}
\pnode(15,9){C}
\pnode(18,8){D}
\ncarc[arrowsize=5pt,arcangleA=45,arcangleB=45,nodesep=5pt]{->}{A}{C}
\ncarc[arrowsize=5pt,arcangleA=-45,arcangleB=-45,nodesep=5pt]{->}{B}{C}
\ncarc[arrowsize=5pt,arcangleA=-45,arcangleB=-45,nodesep=5pt]{->}{A}{B}
\ncarc[arrowsize=5pt,arcangleA=90,arcangleB=90,nodesep=5pt]{->}{A}{D}

\pnode(9,1){E}
\pnode(12,2){F}
\pnode(15,3){G}
\pnode(18,4){H}
\ncarc[arrowsize=5pt,arcangleA=45,arcangleB=45,nodesep=5pt]{->}{E}{G}
\ncarc[arrowsize=5pt,arcangleA=-45,arcangleB=-45,nodesep=5pt]{->}{F}{G}
\ncarc[arrowsize=5pt,arcangleA=-45,arcangleB=-45,nodesep=5pt]{->}{E}{F}
\ncarc[arrowsize=5pt,arcangleA=90,arcangleB=90,nodesep=5pt]{->}{E}{H}

	\psline[linewidth=0.1pt](0,2)(2,0)\psline[linewidth=0.1pt](0,4)(4,0)\psline[linewidth=0.1pt](0,6)(6,0)\psline[linewidth=0.1pt](0,8)(8,0)\psline[linewidth=0.1pt](0,10)(10,0)\psline[linewidth=0.1pt](0,12)(12,0)\psline[linewidth=0.1pt](2,12)(14,0)\psline[linewidth=0.1pt](4,12)(16,0)\psline[linewidth=0.1pt](6,12)(18,0)\psline[linewidth=0.1pt](8,12)(20,0)\psline[linewidth=0.1pt](10,12)(22,0)\psline[linewidth=0.1pt](12,12)(24,0)\psline[linewidth=0.1pt](14,12)(24,2)\psline[linewidth=0.1pt](16,12)(24,4)\psline[linewidth=0.1pt](18,12)(24,6)\psline[linewidth=0.1pt](20,12)(24,8)\psline[linewidth=0.1pt](22,12)(24,10)\psline[linewidth=0.1pt](24,12)(24,12)\psline[linewidth=0.1pt](0,10)(2,12)\psline[linewidth=0.1pt](0,8)(4,12)\psline[linewidth=0.1pt](0,6)(6,12)\psline[linewidth=0.1pt](0,4)(8,12)\psline[linewidth=0.1pt](0,2)(10,12)\psline[linewidth=0.1pt](0,0)(12,12)\psline[linewidth=0.1pt](2,0)(14,12)\psline[linewidth=0.1pt](4,0)(16,12)\psline[linewidth=0.1pt](6,0)(18,12)\psline[linewidth=0.1pt](8,0)(20,12)\psline[linewidth=0.1pt](10,0)(22,12)\psline[linewidth=0.1pt](12,0)(24,12)\psline[linewidth=0.1pt](14,0)(24,10)\psline[linewidth=0.1pt](16,0)(24,8)\psline[linewidth=0.1pt](18,0)(24,6)\psline[linewidth=0.1pt](20,0)(24,4)\psline[linewidth=0.1pt](22,0)(24,2)\psline[linewidth=0.1pt](24,0)(24,0)\psline[linewidth=0.1pt](0,0)(24,0)\psline[linewidth=0.1pt](0,1)(24,1)\psline[linewidth=0.1pt](0,2)(24,2)\psline[linewidth=0.1pt](0,3)(24,3)\psline[linewidth=0.1pt](0,4)(24,4)\psline[linewidth=0.1pt](0,5)(24,5)\psline[linewidth=0.1pt](0,6)(24,6)\psline[linewidth=0.1pt](0,7)(24,7)\psline[linewidth=0.1pt](0,8)(24,8)\psline[linewidth=0.1pt](0,9)(24,9)\psline[linewidth=0.1pt](0,10)(24,10)\psline[linewidth=0.1pt](0,11)(24,11)\psline[linewidth=0.1pt](0,12)(24,12)\psdots[dotstyle=*](0,0) \uput[330](0,0){$(9,3)$}\psdots[dotstyle=*](0,2) \uput[330](0,2){$(8,4)$}\psdots[dotstyle=*](0,4) \uput[330](0,4){$(7,5)$}\psdots[dotstyle=*](0,6) \uput[330](0,6){$(6,6)$}\psdots[dotstyle=*](0,8) \uput[330](0,8){$(5,7)$}\psdots[dotstyle=*](0,10) \uput[330](0,10){$(4,8)$}\psdots[dotstyle=*](0,12) \uput[330](0,12){$(3,9)$}\psdots[dotstyle=*](3,1) \uput[330](3,1){$(7,2)$}\psdots[dotstyle=*](3,3) \uput[330](3,3){$(6,3)$}\psdots[dotstyle=*](3,5) \uput[330](3,5){$(5,4)$}\psdots[dotstyle=*](3,7) \uput[330](3,7){$(4,5)$}\psdots[dotstyle=*](3,9) \uput[330](3,9){$(3,6)$}\psdots[dotstyle=*](3,11) \uput[330](3,11){$(2,7)$}\psdots[dotstyle=*](6,0) \uput[330](6,0){$(6,0)$}\psdots[dotstyle=*](6,2) \uput[330](6,2){$(5,1)$}\psdots[dotstyle=*](6,4) \uput[330](6,4){$(4,2)$}\psdots[dotstyle=*](6,6) \uput[330](6,6){$(3,3)$}\psdots[dotstyle=*](6,8) \uput[330](6,8){$(2,4)$}\psdots[dotstyle=*](6,10) \uput[330](6,10){$(1,5)$}\psdots[dotstyle=*](6,12) \uput[330](6,12){$(0,6)$}\psdots[dotstyle=*](9,1) \uput[330](9,1){$(4,-1)$}\psdots[dotstyle=*](9,3) \uput[330](9,3){$(3,0)$}\psdots[dotstyle=*](9,5) \uput[330](9,5){$(2,1)$}\psdots[dotstyle=*](9,7) \uput[330](9,7){$(1,2)$}\psdots[dotstyle=*](9,9) \uput[330](9,9){$(0,3)$}\psdots[dotstyle=*](9,11) \uput[330](9,11){$(-1,4)$}\psdots[dotstyle=*](12,0) \uput[330](12,0){$(3,-3)$}\psdots[dotstyle=*](12,2) \uput[330](12,2){$(2,-2)$}\psdots[dotstyle=*](12,4) \uput[330](12,4){$(1,-1)$}\psdots[dotstyle=*](12,6) \uput[330](12,6){$(0,0)$}\psdots[dotstyle=*](12,8) \uput[330](12,8){$(-1,1)$}\psdots[dotstyle=*](12,10) \uput[330](12,10){$(-2,2)$}\psdots[dotstyle=*](12,12) \uput[330](12,12){$(-3,3)$}\psdots[dotstyle=*](15,1) \uput[330](15,1){$(1,-4)$}\psdots[dotstyle=*](15,3) \uput[330](15,3){$(0,-3)$}\psdots[dotstyle=*](15,5) \uput[330](15,5){$(-1,-2)$}\psdots[dotstyle=*](15,7) \uput[330](15,7){$(-2,-1)$}\psdots[dotstyle=*](15,9) \uput[330](15,9){$(-3,0)$}\psdots[dotstyle=*](15,11) \uput[330](15,11){$(-4,1)$}\psdots[dotstyle=*](18,0) \uput[330](18,0){$(0,-6)$}\psdots[dotstyle=*](18,2) \uput[330](18,2){$(-1,-5)$}\psdots[dotstyle=*](18,4) \uput[330](18,4){$(-2,-4)$}\psdots[dotstyle=*](18,6) \uput[330](18,6){$(-3,-3)$}\psdots[dotstyle=*](18,8) \uput[330](18,8){$(-4,-2)$}\psdots[dotstyle=*](18,10) \uput[330](18,10){$(-5,-1)$}\psdots[dotstyle=*](18,12) \uput[330](18,12){$(-6,0)$}\psdots[dotstyle=*](21,1) \uput[330](21,1){$(-2,-7)$}\psdots[dotstyle=*](21,3) \uput[330](21,3){$(-3,-6)$}\psdots[dotstyle=*](21,5) \uput[330](21,5){$(-4,-5)$}\psdots[dotstyle=*](21,7) \uput[330](21,7){$(-5,-4)$}\psdots[dotstyle=*](21,9) \uput[330](21,9){$(-6,-3)$}\psdots[dotstyle=*](21,11) \uput[330](21,11){$(-7,-2)$}\psdots[dotstyle=*](24,0) \uput[330](24,0){$(-3,-9)$}\psdots[dotstyle=*](24,2) \uput[330](24,2){$(-4,-8)$}\psdots[dotstyle=*](24,4) \uput[330](24,4){$(-5,-7)$}\psdots[dotstyle=*](24,6) \uput[330](24,6){$(-6,-6)$}\psdots[dotstyle=*](24,8) \uput[330](24,8){$(-7,-5)$}\psdots[dotstyle=*](24,10) \uput[330](24,10){$(-8,-4)$}\psdots[dotstyle=*](24,12) \uput[330](24,12){$(-9,-3)$} 
	\end{pspicture}
    \end{center}
\caption{The movement of non-stationary points in $\Delta_{5}$ for $a=3$ }\label{movement}
\end{figure}

The set of points formed by the union of the affine spaces in the $a$-paving that are
indexed by the vertices in $\Delta_{i}$
is equal to the set of points in the union of $I$-orbits of vertices in $\Delta_{i}$.
To see this, we first refer to Lemma \ref{123 movement}.  In that lemma, we prove that 
some non-stationary points originate from the $I$-orbit of a type $1^{0}$ vertex 
and move to the $I^{a}$-orbit of a type $1^{a}$, $2^{a}$, or $3^{a}$ vertex.  Relative 
to $\Delta_{i}$, these non-stationary points move along the upper edge of $\Delta_{i}$
excluding the vertices of the triangle.  (See figure \ref{movement}.)
In Lemma \ref{567 movement}, we prove the other possibility is that 
non-stationary points originate from the $I$-orbit of a type $7^{0}$ vertex 
and move to the $I^{a}$-orbit of a type $5^{a}$, $6^{a}$, or $7^{a}$ vertex.
 Relative 
to $\Delta_{i}$, these non-stationary points move along the lower edge of $\Delta_{i}$
again excluding the vertices of the triangle.    

We will now show that if the union of the affine spaces indexed by the vertices 
$v_{0}, v_{1}, \ldots, v_{l}$ in
$\Delta_{i-1}$ is closed then there is at least one way to order the 
vertices $v_{l+1}, v_{l+2}, \ldots, v_{k}$ on the boundary of $\Delta_{i}$, such that
\[ \mathbb{A}_{0} \subset \mathbb{A}_{0} \cup \mathbb{A}_{1} \subset \dots \subset \bigcup_{j=0}^{k} \mathbb{A}_{j} \]
are all closed.
Since $\Delta_{0}$ only contains $v_{0}$ and we know that $\mathbb{A}_{0}$ is closed, we 
can proceed by induction and order all the vertices of $X$.

\begin{figure}[!htb]
(i)    \psset{xunit=0.125 in}
    \psset{yunit=0.21650635094611 in}
    \begin{pspicture}[0.9](0,0)(16,8)
    \scriptsize	
	\pspolygon[fillstyle=solid,fillcolor=gray,linestyle=none](6,4)(5,5)(4,4)
	\psdots[dotstyle=*](6,4) \uput[330](6,4){$(0,0)$}
	\psline[linewidth=0.8pt,linestyle=dashed](2,0)(2,8)(14,4)(2,0)
	\psline[linewidth=1.8pt,arrows=o-o](2,0)(2,8)
\psline[linewidth=0.1pt](0,2)(2,0)\psline[linewidth=0.1pt](0,4)(4,0)\psline[linewidth=0.1pt](0,6)(6,0)\psline[linewidth=0.1pt](0,8)(8,0)\psline[linewidth=0.1pt](2,8)(10,0)\psline[linewidth=0.1pt](4,8)(12,0)\psline[linewidth=0.1pt](6,8)(14,0)\psline[linewidth=0.1pt](8,8)(16,0)\psline[linewidth=0.1pt](10,8)(16,2)\psline[linewidth=0.1pt](12,8)(16,4)\psline[linewidth=0.1pt](14,8)(16,6)\psline[linewidth=0.1pt](16,8)(16,8)\psline[linewidth=0.1pt](0,6)(2,8)\psline[linewidth=0.1pt](0,4)(4,8)\psline[linewidth=0.1pt](0,2)(6,8)\psline[linewidth=0.1pt](0,0)(8,8)\psline[linewidth=0.1pt](2,0)(10,8)\psline[linewidth=0.1pt](4,0)(12,8)\psline[linewidth=0.1pt](6,0)(14,8)\psline[linewidth=0.1pt](8,0)(16,8)\psline[linewidth=0.1pt](10,0)(16,6)\psline[linewidth=0.1pt](12,0)(16,4)\psline[linewidth=0.1pt](14,0)(16,2)\psline[linewidth=0.1pt](16,0)(16,0)\psline[linewidth=0.1pt](0,0)(16,0)\psline[linewidth=0.1pt](0,1)(16,1)\psline[linewidth=0.1pt](0,2)(16,2)\psline[linewidth=0.1pt](0,3)(16,3)\psline[linewidth=0.1pt](0,4)(16,4)\psline[linewidth=0.1pt](0,5)(16,5)\psline[linewidth=0.1pt](0,6)(16,6)\psline[linewidth=0.1pt](0,7)(16,7)\psline[linewidth=0.1pt](0,8)(16,8) 
	\end{pspicture}
(ii)    \psset{xunit=0.125 in}
    \psset{yunit=0.21650635094611 in}
    \begin{pspicture}[0.9](0,0)(16,8)
    \scriptsize	
	\pspolygon[fillstyle=solid,fillcolor=gray,linestyle=none](6,4)(5,5)(4,4)
	\psdots[dotstyle=*](6,4) \uput[330](6,4){$(0,0)$}

	\psline[linewidth=0.8pt,linestyle=dashed](2,0)(2,8)(14,4)(2,0)
	\psline[linewidth=1.8pt,arrows=o-o](2,0)(14,4)
	\psline[linewidth=1.8pt,arrows=o-o](2,8)(14,4)
\psline[linewidth=0.1pt](0,2)(2,0)\psline[linewidth=0.1pt](0,4)(4,0)\psline[linewidth=0.1pt](0,6)(6,0)\psline[linewidth=0.1pt](0,8)(8,0)\psline[linewidth=0.1pt](2,8)(10,0)\psline[linewidth=0.1pt](4,8)(12,0)\psline[linewidth=0.1pt](6,8)(14,0)\psline[linewidth=0.1pt](8,8)(16,0)\psline[linewidth=0.1pt](10,8)(16,2)\psline[linewidth=0.1pt](12,8)(16,4)\psline[linewidth=0.1pt](14,8)(16,6)\psline[linewidth=0.1pt](16,8)(16,8)\psline[linewidth=0.1pt](0,6)(2,8)\psline[linewidth=0.1pt](0,4)(4,8)\psline[linewidth=0.1pt](0,2)(6,8)\psline[linewidth=0.1pt](0,0)(8,8)\psline[linewidth=0.1pt](2,0)(10,8)\psline[linewidth=0.1pt](4,0)(12,8)\psline[linewidth=0.1pt](6,0)(14,8)\psline[linewidth=0.1pt](8,0)(16,8)\psline[linewidth=0.1pt](10,0)(16,6)\psline[linewidth=0.1pt](12,0)(16,4)\psline[linewidth=0.1pt](14,0)(16,2)\psline[linewidth=0.1pt](16,0)(16,0)\psline[linewidth=0.1pt](0,0)(16,0)\psline[linewidth=0.1pt](0,1)(16,1)\psline[linewidth=0.1pt](0,2)(16,2)\psline[linewidth=0.1pt](0,3)(16,3)\psline[linewidth=0.1pt](0,4)(16,4)\psline[linewidth=0.1pt](0,5)(16,5)\psline[linewidth=0.1pt](0,6)(16,6)\psline[linewidth=0.1pt](0,7)(16,7)\psline[linewidth=0.1pt](0,8)(16,8) 
	\end{pspicture}
(iii)    \psset{xunit=0.125 in}
    \psset{yunit=0.21650635094611 in}
    \begin{pspicture}[0.9](0,0)(16,8)
    \scriptsize	

	\pspolygon[fillstyle=solid,fillcolor=gray,linestyle=none](6,4)(5,5)(4,4)
	\psdots[dotstyle=*](6,4) \uput[330](6,4){$(0,0)$}

	\psline[linewidth=0.8pt,linestyle=dashed](2,0)(2,8)(14,4)(2,0)
	\psdots[dotstyle=*,dotsize=6pt](2,0)(2,8)(14,4)
\psline[linewidth=0.1pt](0,2)(2,0)\psline[linewidth=0.1pt](0,4)(4,0)\psline[linewidth=0.1pt](0,6)(6,0)\psline[linewidth=0.1pt](0,8)(8,0)\psline[linewidth=0.1pt](2,8)(10,0)\psline[linewidth=0.1pt](4,8)(12,0)\psline[linewidth=0.1pt](6,8)(14,0)\psline[linewidth=0.1pt](8,8)(16,0)\psline[linewidth=0.1pt](10,8)(16,2)\psline[linewidth=0.1pt](12,8)(16,4)\psline[linewidth=0.1pt](14,8)(16,6)\psline[linewidth=0.1pt](16,8)(16,8)\psline[linewidth=0.1pt](0,6)(2,8)\psline[linewidth=0.1pt](0,4)(4,8)\psline[linewidth=0.1pt](0,2)(6,8)\psline[linewidth=0.1pt](0,0)(8,8)\psline[linewidth=0.1pt](2,0)(10,8)\psline[linewidth=0.1pt](4,0)(12,8)\psline[linewidth=0.1pt](6,0)(14,8)\psline[linewidth=0.1pt](8,0)(16,8)\psline[linewidth=0.1pt](10,0)(16,6)\psline[linewidth=0.1pt](12,0)(16,4)\psline[linewidth=0.1pt](14,0)(16,2)\psline[linewidth=0.1pt](16,0)(16,0)\psline[linewidth=0.1pt](0,0)(16,0)\psline[linewidth=0.1pt](0,1)(16,1)\psline[linewidth=0.1pt](0,2)(16,2)\psline[linewidth=0.1pt](0,3)(16,3)\psline[linewidth=0.1pt](0,4)(16,4)\psline[linewidth=0.1pt](0,5)(16,5)\psline[linewidth=0.1pt](0,6)(16,6)\psline[linewidth=0.1pt](0,7)(16,7)\psline[linewidth=0.1pt](0,8)(16,8) 
	\end{pspicture}
\caption{Three stages in ordering the vertices of $\Delta_{i}$}\label{order steps}
\end{figure}
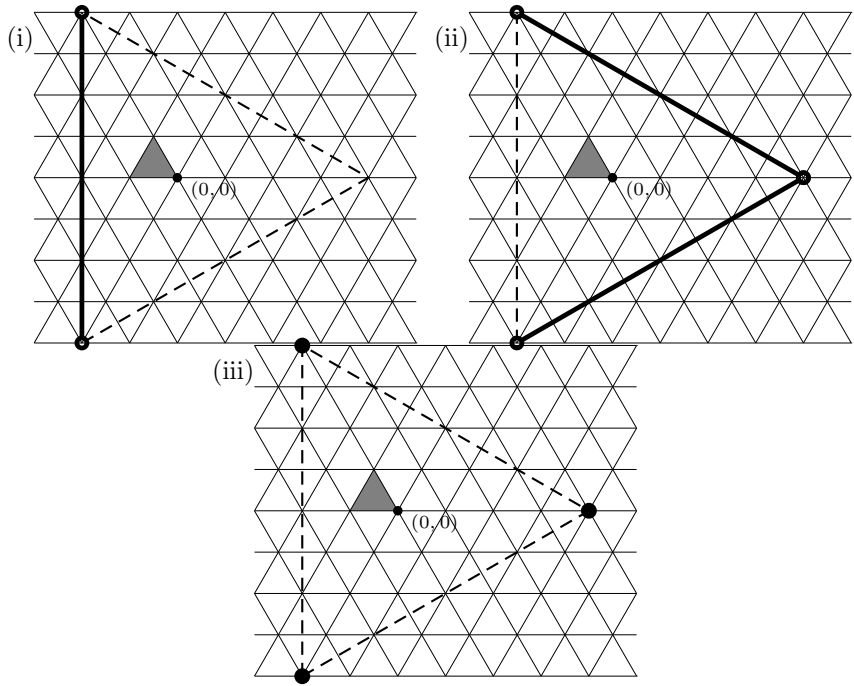

Assume the union of the affine spaces indexed by the vertices in
$\Delta_{i-1}$ is closed.  We will describe a valid order 
for the vertices on the boundary of $\Delta_{i}$ in three stages.
Figure \ref{order steps} graphically depicts the three stages.

In stage (i), we add the vertices along the vertical edge of $\Delta_{i}$ excluding
the two vertices at the ends of this edge.  These vertices can be ordered
so that the dimensions of the associated affine spaces forms a weakly increasing sequence.
To validate the order, we must show that for each $v_{i}$ in stage (i) the set
\[ \bigcup_{j=0}^{i} \mathbb{A}_{j} \]
is closed.  We will use the geometry of minimal galleries to show the union is closed.

Suppose that $v_{i}$ is a vertex from stage (i).  Take a minimal gallery $M$ connecting
$v$ to the base alcove (the alcove stabilized by $I$).  We can express $M$ as a word
of simple reflections $M = s_{1} \ldots s_{j}$.  Each subword of $M$ represents
a vertex in $X$.  If every subword of $M$ represents a vertex that is in
\[ \bigcup_{j=0}^{i-1} \mathbb{A}_{j} \]
then
\[ \bigcup_{j=0}^{i} \mathbb{A}_{j} \] is closed.

\begin{figure}[!htb]
    \psset{xunit=0.125 in}
    \psset{yunit=0.21650635094611 in}
    \begin{pspicture}[0.9](0,0)(16,8)
    \scriptsize	

\psdots[dotstyle=*](2,6)	\pspolygon[fillstyle=solid,fillcolor=gray,linestyle=none](6,4)
(4,4)(2,6)(4,6)
	\psdots[dotstyle=*](6,4) \uput[330](6,4){$(0,0)$}

	\psline[linewidth=0.8pt,linestyle=dashed](2,0)(2,8)(14,4)(2,0)
	\psline[linewidth=0.8pt,linestyle=dashed](3,1)(3,7)(12,4)(3,1)

\psline[linewidth=0.1pt](0,2)(2,0)\psline[linewidth=0.1pt](0,4)(4,0)\psline[linewidth=0.1pt](0,6)(6,0)\psline[linewidth=0.1pt](0,8)(8,0)\psline[linewidth=0.1pt](2,8)(10,0)\psline[linewidth=0.1pt](4,8)(12,0)\psline[linewidth=0.1pt](6,8)(14,0)\psline[linewidth=0.1pt](8,8)(16,0)\psline[linewidth=0.1pt](10,8)(16,2)\psline[linewidth=0.1pt](12,8)(16,4)\psline[linewidth=0.1pt](14,8)(16,6)\psline[linewidth=0.1pt](16,8)(16,8)\psline[linewidth=0.1pt](0,6)(2,8)\psline[linewidth=0.1pt](0,4)(4,8)\psline[linewidth=0.1pt](0,2)(6,8)\psline[linewidth=0.1pt](0,0)(8,8)\psline[linewidth=0.1pt](2,0)(10,8)\psline[linewidth=0.1pt](4,0)(12,8)\psline[linewidth=0.1pt](6,0)(14,8)\psline[linewidth=0.1pt](8,0)(16,8)\psline[linewidth=0.1pt](10,0)(16,6)\psline[linewidth=0.1pt](12,0)(16,4)\psline[linewidth=0.1pt](14,0)(16,2)\psline[linewidth=0.1pt](16,0)(16,0)\psline[linewidth=0.1pt](0,0)(16,0)\psline[linewidth=0.1pt](0,1)(16,1)\psline[linewidth=0.1pt](0,2)(16,2)\psline[linewidth=0.1pt](0,3)(16,3)\psline[linewidth=0.1pt](0,4)(16,4)\psline[linewidth=0.1pt](0,5)(16,5)\psline[linewidth=0.1pt](0,6)(16,6)\psline[linewidth=0.1pt](0,7)(16,7)\psline[linewidth=0.1pt](0,8)(16,8) 
	\end{pspicture}
\caption{Example of a minimal gallery for stage (i)}\label{minimal gallery}
\end{figure}
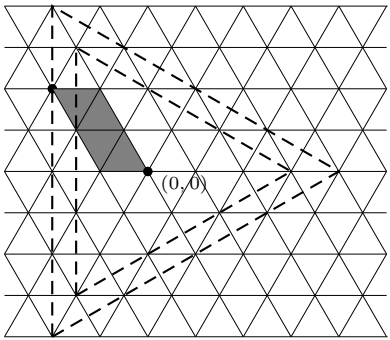

Since each affine space associated to a vertex in stage (i) is simply the $I$-orbit
of that vertex, the dimension of the affine space is the length of the minimal gallery 
connecting the vertex to the alcove stabilized by $I$.  Thus, we are proposing that
the vertices in stage (i) can ordered by the length of the minimal gallery associated 
to the vertex.  By considering the possible subwords of a minimal gallery associated to
a vertex in stage (i), we see that the vertex represented by the subword is either
another stage (i) vertex or a vertex in $\Delta_{i-1}$.  See figure \ref{minimal gallery}.
Since we order the vertices in stage (i) by the length of the associated minimal gallery,
the proposed order is valid.

In stage (ii), we can use the dimension of the $I^{a}$-orbit of each 
vertex to order the vertices.  To prove this, it is enough to show that the union $Y$
of the affine spaces associated to $\Delta_{i-1}$, stage (i), and stage (ii) 
is closed.   Then, since the union
of the affine spaces associated to $\Delta_{i-1}$ and stage (i)
is closed, the  union
of the affine spaces associated to stage (ii) is open in $Y$.  
Therefore, we can work in the topological subspace 
defined by the union
of the affine spaces associated to stage (ii).  Each of those affine
spaces is a subspace of the $I^{a}$-orbit of the associated vertex. 
Hence, if we order the vertices such that the dimension the $I^{a}$-orbit
of each vertex increases, the order is valid.

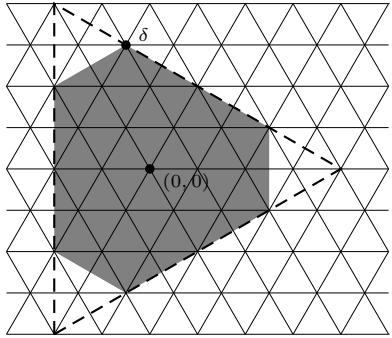
\begin{figure}[!htb]
    \psset{xunit=0.125 in}
    \psset{yunit=0.21650635094611 in}
    \begin{pspicture}[0.9](0,0)(16,8)
    \scriptsize	

\pspolygon[fillstyle=solid,fillcolor=gray,linestyle=none](5,7)(2,6)(2,2)(5,1)(11,3)(11,5)
	\psdots[dotstyle=*](6,4) \uput[330](6,4){$(0,0)$}
	\psdots[dotstyle=*](5,7) \uput[30](5,7){$\delta$}

	\psline[linewidth=0.8pt,linestyle=dashed](2,0)(2,8)(14,4)(2,0)

\psline[linewidth=0.1pt](0,2)(2,0)\psline[linewidth=0.1pt](0,4)(4,0)\psline[linewidth=0.1pt](0,6)(6,0)\psline[linewidth=0.1pt](0,8)(8,0)\psline[linewidth=0.1pt](2,8)(10,0)\psline[linewidth=0.1pt](4,8)(12,0)\psline[linewidth=0.1pt](6,8)(14,0)\psline[linewidth=0.1pt](8,8)(16,0)\psline[linewidth=0.1pt](10,8)(16,2)\psline[linewidth=0.1pt](12,8)(16,4)\psline[linewidth=0.1pt](14,8)(16,6)\psline[linewidth=0.1pt](16,8)(16,8)\psline[linewidth=0.1pt](0,6)(2,8)\psline[linewidth=0.1pt](0,4)(4,8)\psline[linewidth=0.1pt](0,2)(6,8)\psline[linewidth=0.1pt](0,0)(8,8)\psline[linewidth=0.1pt](2,0)(10,8)\psline[linewidth=0.1pt](4,0)(12,8)\psline[linewidth=0.1pt](6,0)(14,8)\psline[linewidth=0.1pt](8,0)(16,8)\psline[linewidth=0.1pt](10,0)(16,6)\psline[linewidth=0.1pt](12,0)(16,4)\psline[linewidth=0.1pt](14,0)(16,2)\psline[linewidth=0.1pt](16,0)(16,0)\psline[linewidth=0.1pt](0,0)(16,0)\psline[linewidth=0.1pt](0,1)(16,1)\psline[linewidth=0.1pt](0,2)(16,2)\psline[linewidth=0.1pt](0,3)(16,3)\psline[linewidth=0.1pt](0,4)(16,4)\psline[linewidth=0.1pt](0,5)(16,5)\psline[linewidth=0.1pt](0,6)(16,6)\psline[linewidth=0.1pt](0,7)(16,7)\psline[linewidth=0.1pt](0,8)(16,8) 
	\end{pspicture}
\caption{The convex hull of the Weyl group orbit of $\delta$}\label{w-orbit}
\end{figure}

To prove that the union
of the affine spaces associated to $\Delta_{i-1}$, stage (i), and stage (ii) 
is closed, we use the fact that the union of the affine spaces associated to
the vertices in the convex hull of the Weyl group orbit of a vertex $\delta$
is closed.  If we let $\delta$ be the vertex in stage (ii) that is adjacent
to the top most geometric vertex of $\Delta_{i}$, then from figure \ref{w-orbit}
we can see the convex hull of the Weyl group orbit of $\delta$ is the 
set of vertices in  $\Delta_{i-1}$, stage (i), and stage (ii). 

Stage (iii) is the three geometric vertices of $\Delta_{i}$.  Since 
the union of the affine spaces indexed by the vertices in
$\Delta_{i}$ is identical to the union on the $I$-orbits of those vertices and since
the affine spaces associated to the vertices from stage (iii) are $I$-orbits 
of those vertices, we can simply order those vertices such that the 
dimensions of the associated affine spaces forms a weakly increasing sequence.

\begin{figure}[!htb]
	\begin{center}
    \psset{xunit=0.08588085254195725 in}
    \psset{yunit=0.14875 in}
    \begin{pspicture}(0,0)(30,20)
    \tiny

\pspolygon[fillstyle=solid,fillcolor=gray,linestyle=none](12,10)(11,11)(10,10)
\pspolygon[fillstyle=solid,fillcolor=gray,linestyle=none](4,10)(3,11)(2,10)

	\psline(3,19)(3,1)(30,10)(3,19)

{\small
\uput[180](3,11){\fbox{1}}
\uput[180](3,9){\fbox{2}}
\uput[180](3,13){\fbox{3}}
\uput[180](3,7){\fbox{4}}
\uput[180](3,15){\fbox{5}}
\uput[180](3,5){\fbox{6}}
\uput[180](3,17){\fbox{7}}
\uput[180](3,3){\fbox{8}}

\uput[30](9,17){\fbox{9}}
\uput[30](12,16){\fbox{11}}
\uput[30](6,18){\fbox{13}}
\uput[30](15,15){\fbox{15}}
\uput[30](18,14){\fbox{17}}
\uput[30](21,13){\fbox{19}}
\uput[30](24,12){\fbox{21}}
\uput[30](27,11){\fbox{23}}

\uput[330](9,3){\fbox{10}}
\uput[330](12,4){\fbox{12}}
\uput[330](6,2){\fbox{14}}
\uput[330](15,5){\fbox{16}}
\uput[330](18,6){\fbox{18}}
\uput[330](21,7){\fbox{20}}
\uput[330](24,8){\fbox{22}}
\uput[330](27,9){\fbox{24}}

\uput[180](3,19){\fbox{25}}
\uput[180](3,1){\fbox{26}}
\uput[0](30,10){\fbox{27}}
}

	\psline[linewidth=0.1pt](0,2)(2,0)\psline[linewidth=0.1pt](0,4)(4,0)\psline[linewidth=0.1pt](0,6)(6,0)\psline[linewidth=0.1pt](0,8)(8,0)\psline[linewidth=0.1pt](0,10)(10,0)\psline[linewidth=0.1pt](0,12)(12,0)\psline[linewidth=0.1pt](0,14)(14,0)\psline[linewidth=0.1pt](0,16)(16,0)\psline[linewidth=0.1pt](0,18)(18,0)\psline[linewidth=0.1pt](0,20)(20,0)\psline[linewidth=0.1pt](2,20)(22,0)\psline[linewidth=0.1pt](4,20)(24,0)\psline[linewidth=0.1pt](6,20)(26,0)\psline[linewidth=0.1pt](8,20)(28,0)\psline[linewidth=0.1pt](10,20)(30,0)\psline[linewidth=0.1pt](12,20)(30,2)\psline[linewidth=0.1pt](14,20)(30,4)\psline[linewidth=0.1pt](16,20)(30,6)\psline[linewidth=0.1pt](18,20)(30,8)\psline[linewidth=0.1pt](20,20)(30,10)\psline[linewidth=0.1pt](22,20)(30,12)\psline[linewidth=0.1pt](24,20)(30,14)\psline[linewidth=0.1pt](26,20)(30,16)\psline[linewidth=0.1pt](28,20)(30,18)\psline[linewidth=0.1pt](30,20)(30,20)\psline[linewidth=0.1pt](0,18)(2,20)\psline[linewidth=0.1pt](0,16)(4,20)\psline[linewidth=0.1pt](0,14)(6,20)\psline[linewidth=0.1pt](0,12)(8,20)\psline[linewidth=0.1pt](0,10)(10,20)\psline[linewidth=0.1pt](0,8)(12,20)\psline[linewidth=0.1pt](0,6)(14,20)\psline[linewidth=0.1pt](0,4)(16,20)\psline[linewidth=0.1pt](0,2)(18,20)\psline[linewidth=0.1pt](0,0)(20,20)\psline[linewidth=0.1pt](2,0)(22,20)\psline[linewidth=0.1pt](4,0)(24,20)\psline[linewidth=0.1pt](6,0)(26,20)\psline[linewidth=0.1pt](8,0)(28,20)\psline[linewidth=0.1pt](10,0)(30,20)\psline[linewidth=0.1pt](12,0)(30,18)\psline[linewidth=0.1pt](14,0)(30,16)\psline[linewidth=0.1pt](16,0)(30,14)\psline[linewidth=0.1pt](18,0)(30,12)\psline[linewidth=0.1pt](20,0)(30,10)\psline[linewidth=0.1pt](22,0)(30,8)\psline[linewidth=0.1pt](24,0)(30,6)\psline[linewidth=0.1pt](26,0)(30,4)\psline[linewidth=0.1pt](28,0)(30,2)\psline[linewidth=0.1pt](30,0)(30,0)\psline[linewidth=0.1pt](0,0)(30,0)\psline[linewidth=0.1pt](0,1)(30,1)\psline[linewidth=0.1pt](0,2)(30,2)\psline[linewidth=0.1pt](0,3)(30,3)\psline[linewidth=0.1pt](0,4)(30,4)\psline[linewidth=0.1pt](0,5)(30,5)\psline[linewidth=0.1pt](0,6)(30,6)\psline[linewidth=0.1pt](0,7)(30,7)\psline[linewidth=0.1pt](0,8)(30,8)\psline[linewidth=0.1pt](0,9)(30,9)\psline[linewidth=0.1pt](0,10)(30,10)\psline[linewidth=0.1pt](0,11)(30,11)\psline[linewidth=0.1pt](0,12)(30,12)\psline[linewidth=0.1pt](0,13)(30,13)\psline[linewidth=0.1pt](0,14)(30,14)\psline[linewidth=0.1pt](0,15)(30,15)\psline[linewidth=0.1pt](0,16)(30,16)\psline[linewidth=0.1pt](0,17)(30,17)\psline[linewidth=0.1pt](0,18)(30,18)\psline[linewidth=0.1pt](0,19)(30,19)\psline[linewidth=0.1pt](0,20)(30,20)

\psdots[dotstyle=*](12,10) \uput[330](12,10){$(0,0)$}
\psdots[dotstyle=*](4,10) \uput[330](4,10){$(4,4)$}
\psdots[dotstyle=*](3,1)
\psdots[dotstyle=*](3,3)
\psdots[dotstyle=*](3,5)
\psdots[dotstyle=*](3,7)
\psdots[dotstyle=*](3,9)
\psdots[dotstyle=*](3,11)
\psdots[dotstyle=*](3,13)
\psdots[dotstyle=*](3,15)
\psdots[dotstyle=*](3,17)
\psdots[dotstyle=*](3,19)
\psdots[dotstyle=*](6,2)
\psdots[dotstyle=*](6,18)
\psdots[dotstyle=*](9,3)
\psdots[dotstyle=*](9,17)
\psdots[dotstyle=*](12,4)
\psdots[dotstyle=*](12,16)
\psdots[dotstyle=*](15,5)
\psdots[dotstyle=*](15,15)
\psdots[dotstyle=*](18,6)
\psdots[dotstyle=*](18,14)
\psdots[dotstyle=*](21,7)
\psdots[dotstyle=*](21,13)
\psdots[dotstyle=*](24,8)
\psdots[dotstyle=*](24,12)
\psdots[dotstyle=*](27,9)
\psdots[dotstyle=*](27,11)
\psdots[dotstyle=*](30,10)

	\end{pspicture}
    \end{center}
\caption{Order on the vertices of $\Delta_{9}$ for $a=4$}\label{order example}
\end{figure}
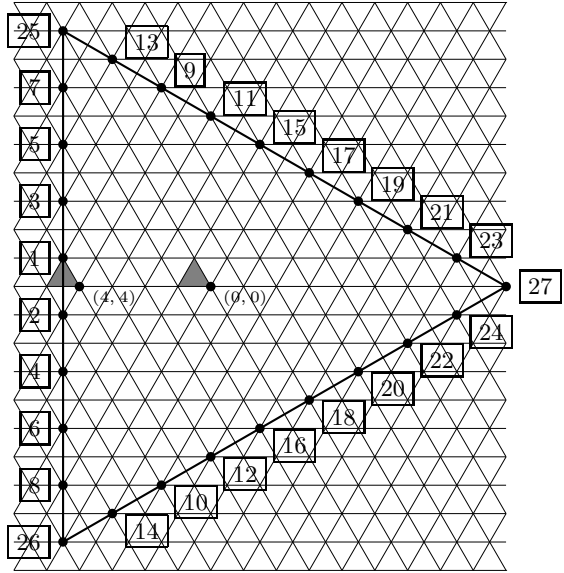

In figure \ref{order example}, the vertices of $\Delta_{9}$ are numbered according
to the order we just described.  In stage (ii), some $I^{a}$-orbits have the same
dimension so figure \ref{order example} gives one of several valid orderings.
Note that the smallest numbers are along the vertical edge,
which corresponds to stage (i), and the three largest numbers label the geometric vertices
of $\Delta_{9}$, which corresponds to stage (iii).
Further, the vertices in stage (ii) are labeled in order of increasing dimension of the 
associated $I^{a}$-orbit of the vertex (which is equal to the length of a minimal gallery between
the vertex and the alcove stabilized by $I^{a}$).  

\begin{thebibliography}{7}
\bibitem[Fan96]{Fan96}
	C. Kenneth Fan, \emph{Euler characteristic of certain affine flag varieties},
	Transform.~Groups {\bf 1} (1996), 35-39.
\bibitem[GKM]{GKM}
	M.~Goresky, R.~Kottwitz and R.~MacPherson
	\emph{Purity of equivalued affine Springer fibers},
	2003, arXiv:math.RT/0305141.
\bibitem[KL88]{KL88}
	D.~Kazhdan and G.~Lusztig,
	\emph{Fixed point varieties on affine flag manifolds},
	Israel J. Math. {\bf 62} (1988), 129-168.
\bibitem[LW]{LW}
	G.~Laumon and J.-L.~Waldspurger, \emph{Sur le lemme fondamental pour les groupes unitaires: 
	le cas totalement ramifi\'e et homog\`ene}, math.AG/9901114.
\bibitem[LS91]{LS91}
	G.~Lusztig and J.~M.~Smelt, \emph{Fixed point varieties on the space of lattices},
	Bull.~London Math.~Soc.~{\bf 23} (1991), 213-218.
\bibitem[Sag00]{Sag00}
	D.~Sage, \emph{The geometry of fixed point varieties on affine flag manifolds},
	Trans.~Amer.~Math.~Soc.~{\bf 352} (2000), 2087-2119.
\bibitem[Som97]{Som97}
	E.~Sommers, \emph{A family of affine Weyl group representations}, Transform.~Groups {\bf 2} (1997), 375-390.
\end{thebibliography}
\end{document}